\documentclass[12pt,oneside,reqno]{amsart}
\usepackage{graphicx}
\usepackage{mathrsfs}
\usepackage{stmaryrd}
\usepackage{amsfonts}
\usepackage{amsfonts,color}
\usepackage{cite}
\usepackage{enumerate,amsmath,amssymb,amsthm}
\usepackage{booktabs} 
\usepackage{diagbox} 

\newcommand{\dif}{\mathrm{d}}

\newcommand{\be}{\begin{eqnarray}}
	\newcommand{\ee}{\end{eqnarray}}
\newcommand{\ce}{\begin{eqnarray*}}
	\newcommand{\de}{\end{eqnarray*}}
\newtheorem{theorem}{Theorem}[section]
\newtheorem{lemma}[theorem]{Lemma}
\newtheorem{remark}[theorem]{Remark}
\newtheorem{definition}[theorem]{Definition}
\newtheorem{proposition}[theorem]{Proposition}
\newtheorem{Examples}[theorem]{Examples}
\newtheorem{corollary}[theorem]{Corollary}
\newtheorem{condition}[theorem]{Condition}
\def\e{\varepsilon}
\def\t{\theta}
\def\a{\alpha}

\def\b{\beta}
\def\d{\delta}
\def\p{\partial}
\def\g{\gamma}

\def\l{\lambda}

\def\[{{\Big[}}
\def\]{{\Big]}}
\def\<{{\langle}}
\def\>{{\rangle}}
\def\({{\Big(}}
\def\){{\Big)}}

\def\no{\nonumber}
\def\bt{\begin{theorem}}
	\def\et{\end{theorem}}
\def\bl{\begin{lemma}}
	\def\el{\end{lemma}}
\def\br{\begin{remark}}
	\def\er{\end{remark}}
\def\bx{\begin{Examples}}
	\def\ex{\end{Examples}}
\def\bd{\begin{definition}}
	\def\ed{\end{definition}}
\def\bp{\begin{proposition}}
	\def\ep{\end{proposition}}
\def\bc{\begin{corollary}}
	\def\ec{\end{corollary}}
\def\bco{\begin{condition}}
	\def\eco{\end{condition}}
\def\cA{{\mathcal A}}

\def\cG{{\mathcal G}}

\def\cK{{\mathcal K}}
\def\cL{{\mathcal L}}
\def\cM{{\mathcal M}}

\def\mE{{\mathbb E}}

\def\mH{{\mathbb H}}

\def\mN{{\mathbb N}}

\def\mP{{\mathbb P}}

\def\mR{{\mathbb R}}

\def\mU{{\mathbb U}}
\def\mV{{\mathbb V}}

\def\mX{{\mathbb X}}

\def\sB{{\mathscr B}}

\def\sF{{\mathscr F}}

\def\sU{{\mathscr U}}

\def\geq{\geqslant}
\def\leq{\leqslant}

\begin{document}
	
\allowdisplaybreaks
\title{Large deviation principles for multiscale stochastic Burgers equations with reflection}
	
\author{Huijie Qiao}

	
	
	
\date{}
	
\dedicatory{School of Mathematics,
		Southeast University\\
		Nanjing, Jiangsu 211189, China\\
		hjqiaogean@seu.edu.cn}
	
\begin{abstract}
This study investigates multiscale stochastic Burgers equations with reflection, wherein the slow component is modeled by a stochastic Burgers equation with reflection and the fast component by a stochastic reaction-diffusion equation with reflection. Using the weak convergence approach, we rigorously establish a large deviation principle for the slow component. Key technical tools include the penalization method, carefully constructed stopping times, and a refined adaptation of Khasminskii's classical time discretization scheme. To illustrate the applicability and significance of our result, we conclude with a concrete example.

\end{abstract}
	
\maketitle \rm
	
\section{Introduction}

Multiscale systems arise naturally across a wide range of scientific disciplines, including physics, chemistry, biology, mathematical finance, and climatology (see, e.g., \cite{ffk, hqs, kk, tll, wtry}). In a multiscale system, one component changes rapidly while the other changes slowly. So it is difficult to control the whole system. Characterizing their asymptotic behavior is therefore essential for both systematic analysis and principled modeling.

The large deviation principle (LDP) constitutes a foundational pillar of asymptotic probability theory, quantifying the exponential decay rates of probabilities associated with rare events. In particular, the Freidlin-Wentzell LDP provides precise estimates for the likelihood that the sample paths of an It\^o diffusion deviate significantly from their deterministic mean trajectories under small-noise perturbations, measured in an appropriate path-space norm (\cite{de}). Freidlin-Wentzell LDPs have been established for a wide range of stochastic equations, including multiscale stochastic differential equations (SDEs) (\cite{abz, ffk, hlls, yK, kp, tll, rL, jll, aap, ks1, aV1, aV2}) and multiscale stochastic partial differential equations (SPDEs) (\cite{ghl, hqs, hlll, hll, hss, mll, swxy, wrd, xxyp, ytw}). 

SDEs with reflection were first introduced by Skorokhod \cite{aS} to characterize diffusion processes confined to a closed domain by means of boundary-localized drift terms. Subsequently, such equations have been widely adopted in applied probability and mathematical modeling, including molecular dynamics, queueing systems, storage models and stochastic control (cf. \cite{btwy, aP}). Motivated by applications to mechanics and stochastic control theory, the framework was extended to infinite-dimensional settings by Haussmann and Pardoux \cite{hp}. Later, Funaki and Olla \cite{fo} studied the wave motion near the hydrodynamic limit of a Ginzburg-Landau $\triangledown\phi$ interface model on a wall, and proved that the fluctuations of a $ \triangledown\phi$ interface model near a hard wall converge in law to the stationary solution of a SPDE with reflection. That is, the evolution of the random interface near the hard wall can be simulated by a SPDE with reflection. Recently, SPDEs with reflection have also emerged as tractable models for order-book dynamics in financial markets (\cite{hkn}). A substantial body of theoretical and analytical results has since been developed for SPDEs with reflection, encompassing existence and uniqueness of solutions, LDPs, invariant measures, and scaling limits (\cite{blz, zbz, my, q, wzz, xz, yz, zhangt0, zhangt}).

For multiscale SDEs with reflection, Kushner \cite{ku} derived a Freidlin-Wentzell LDP via functional occupation measures. Subsequently, in \cite{q0} and \cite{q2}, we addressed multiscale multivalued SDEs to establish corresponding LDPs by the weak convergence approach and the viscosity solution method, respectively. More recently, this analysis was generalized to multiscale multivalued McKean-Vlasov SDEs in \cite{q1}. Notably, when the maximal monotone operator reduces to the subdifferential operator for the indicator function of a closed convex set, the multivalued dynamics collapse precisely to reflected dynamics. Thus, existing results on multiscale SDEs with reflection constitute special cases of the aforementioned frameworks. Nevertheless, the question whether a Freidlin-Wentzell LDP holds for multiscale SPDEs with reflection remains open and is both natural and significant. In this work, we provide an affirmative answer.

Specifically, we consider the following stochastic slow-fast system on the spatial domain $[0,1]$:
\be\left\{\begin{array}{ll}
\frac{\p X_t^{\e,\g}(\xi)}{\p t}=\frac{\p^2 X_t^{\e,\g}(\xi)}{\p\xi^2}+\frac{1}{2}\frac{\p(X_t^{\e,\g})^2(\xi)}{\p \xi}+F_1(X_t^{\e,\g},Y_t^{\e,\g})(\xi)+\sqrt \e G_1(X_t^{\e,\g})\frac{\p W^1(t,\xi)}{\p t}\\
\qquad\qquad\qquad +\frac{\p K^{1,\e,\g}(t,\xi)}{\p t},\\
\frac{\p Y_t^{\e,\g}(\xi)}{\p t}=\frac{1}{\g}\left[\frac{\p^2 Y_t^{\e,\g}(\xi)}{\p\xi^2}+F_2(X_t^{\e,\g},Y_t^{\e,\g})(\xi)\right]+\frac{1}{\sqrt \g} G_2(X_t^{\e,\g},Y_t^{\e,\g})\frac{\p W^2(t,\xi)}{\p t}+\frac{\p K^{2,\e,\g}(t,\xi)}{\p t},\\
X_t^{\e,\g}(\xi)\geq 0, \quad Y_t^{\e,\g}(\xi)\geq 0, \\
X_t^{\e,\g}(0)=X_t^{\e,\g}(1)=Y_t^{\e,\g}(0)=Y_t^{\e,\g}(1)=0, \quad t>0,\\
X_0^{\e,\g}(\xi)=x_0(\xi)\geq 0, \quad Y_0^{\e,\g}(\xi)=y_0(\xi)\geq 0,
\end{array}
\right.
\label{0eq}
\ee
where $\e>0$ denotes a small noise intensity parameter, and $\g=\g(\e)>0$ represents the time-scale separation ratio between the slow component $X^{\e,\g}$ and the fast component $Y^{\e,\g}$. The nonlinearities $F_1$, $G_1$, $F_2$, $G_2$ satisfy regularity conditions and $W^1$, $W^2$ are mutually independent cylindrical Wiener processes defined on a complete filtered probability space $(\Omega, \sF, \{\sF_t\}_{t\geq 0}, \mP)$. To ensure pathwise nonnegativity of the solutions $X^{\e,\g}$ and $Y^{\e,\g}$, reflection random measures $K^{1,\e,\g}, K^{2,\e,\g}$ are added to the system (\ref{0eq}) (specified in Subsection \ref{msbewr}). Since the slow component of system (\ref{0eq}) is governed by a stochastic Burgers equation augmented with boundary reflection at zero, we refer to the full system (\ref{0eq}) as multiscale stochastic Burgers equations with reflection. In the absence of reflection (i.e. without $K^{1,\e,\g}, K^{2,\e,\g}$), both strong and weak averaging principles for this system were established in \cite{dsxz}. Subsequently, Sun et al. \cite{swxy} derived its Freidlin-Wentzell LDP using the weak convergence approach. For the reflected system (\ref{0eq}), Ma and Yang \cite{my} recently proved the strong averaging principle. In this work, we extend that framework to establish the Freidlin-Wentzell LDP for the system (\ref{0eq}) via the weak convergence method.

We emphasize that although the Freidlin-Wentzell LDP and uniform LDP for stochastic Burgers equations with reflection have been established in \cite{wzz} and \cite{q}, respectively, extending these results to multiscale stochastic Burgers equations with reflection is highly nontrivial. First, the averaged skeleton equation (\ref{e0conteq}), which arises in the averaging limit, is structurally more intricate than the skeleton equations considered in \cite{wzz} and \cite{q}. Thus, establishing convergence of its solution demands a refined analytical framework. Second, due to the intrinsic coupling between the slow and fast components in the controlled system (\ref{contequa}), direct verification of the convergence $X^{\varepsilon,\gamma,u_\varepsilon} \to \bar{X}^{u_\varepsilon}$ in $C([0,T], \mH)$ is not feasible. Instead, we develop a robust approximation scheme involving energy estimates and an auxiliary process. Moreover, in contrast to the reflected-free setting analyzed in \cite{swxy}, the presence of the reflection term introduces two fundamental technical obstacles: $(i)$ controlling the $\mH$-norm difference $\sup\limits_{t\in[0,T]}|\bar{X}_t^{h_\varepsilon} - \bar{X}_t^{h}|_{\mH}$ under the weak convergence $h_\varepsilon \to h$ in $L^2([0,T], \sU)$, where the convergence of $\bar{K}^{1,h_\varepsilon}$ to $\bar{K}^{1,h}$ cannot obtained; and $(ii)$ deriving a uniform temporal modulus of continuity estimate for $X^{\varepsilon,\gamma,u_\varepsilon}$, specifically, bounding $|X^{\varepsilon,\gamma,u_\varepsilon}_{t+r} - X^{\varepsilon,\gamma,u_\varepsilon}_t|_{\mH}$ for arbitrary $t > 0$ and arbitrarily small $r>0$, which is complicated by the singular nature of the reflection and its interaction with the multiscale dynamics.

Finally, we outline the main steps of the proof. To verify Condition \ref{cond}$(i)$, we introduce penalization approximations $\bar{X}^{h_\varepsilon,n}$ and $\bar{X}^{h,n}$ of $\bar{X}^{h_\varepsilon}$ and $\bar{X}^{h}$, respectively. We then establish two uniform convergence properties: $(a)$ $\bar{X}^{h,n} \to \bar{X}^{h}$ in $C([0,T],\mH) \cap L^2([0,T],\mV)$, uniformly over all $h \in \mathbf{D}^N$, as $n \to \infty$; and $(b)$ $\bar{X}^{h_\varepsilon,n} \to \bar{X}^{h,n}$ in $C([0,T], \mH)$, as $\varepsilon \to 0$. These convergence results, which are derived via energy estimates, compactness arguments, and stability analysis under the penalization parameter, complete the verification of Condition \ref{cond} $(i)$ with the help of three Sobolev embedding theorems in a more streamlined manner than in prior works \cite{wzz,swxy}. For Condition \ref{cond} $(ii)$, the intrinsic nonlinearity in Eq.(\ref{contequa}) and (\ref{e0conteq}) precludes direct application of standard moment bounds for high-order moments. To overcome this limitation, we introduce two carefully designed stopping times that localize the dynamics to bounded energy regimes while preserving uniform integrability of the associated stochastic integrals. Finally, we adapt Khasminskii's classical time discretization method to control the error term for $\mathbb{E}\sup\limits_{t\in[0,T\wedge \tau_2]} |X^{\varepsilon,\gamma,u_\varepsilon}_t - \bar{X}^{u_\varepsilon}_t|_{\mH}^2$ where $\tau_2$ is defined in (\ref{tau2defi}), leveraging the Gronwall inequality.

The remainder of this paper is structured as follows. Section \ref{pre} collects notations and some results related with LDPs. The precise formulation of our main result appears in Section \ref{main}. Section \ref{fldpthproo} contains the full proof, including technical estimates and convergence arguments. In Section \ref{exam} we give an example to explain our result. Section \ref{app} compiles auxiliary properties of the operators $B$ and $b$ used throughout the analysis.

The following convention will be used throughout the paper: $C$ with or without indices will denote different positive constants whose values may change from one place to another.

\section{Preliminaries}\label{pre}

In this section, we introduce some notations and LDPs.

\subsection{Notation}

In this subsection, we introduce some notations used in the sequel.

Let $C_c([0,1])$ be the set of all continuous functions on $[0,1]$ with compact supports. Let $C^2_c([0,1])$ and  $C^\infty_c([0,1])$  be the subsets of $C_c([0,1])$ where all functions have $2$ and infinite order derivatives, respectively. 

Let $\mH=L^2([0,1], \mR)$ be the usual $L^2$-space with the norm $|\cdot|_\mH$ and inner product $\<\cdot, \cdot\>_\mH$. Denote by $\mV$ the Sobolev space of order one, i.e., $\mV$ is the completion of $C_c^{\infty}([0,1])$ under the norm $\|u\|_\mV^2=\int_0^1\left(\frac{\partial u}{\partial x}\right)^2 d x$. Remark that $\mV=\left\{u \in \mH^{1}([0,1]): u(0)=u(1)=0\right\}$, where $\mH^{1}([0,1])$ denotes the usual Sobolev space of absolutely continuous functions defined on $[0,1]$ whose derivatives belong to $\mH$. $\mV^*$ denotes the dual space of $\mV$ and the dualization between $\mV$ and $\mV^*$ is denoted by ${_{\mV^*}}\<\cdot,\cdot\>_\mV$. 

Let $\mU_i$ be the Hilbert space with norm $|\cdot|_{\mU_i}$ and inner product $\<\cdot, \cdot\>_{\mU_i}, i=1,2$. Let $\cL_2(\mU_i,\mH)$ be the collection of all Hilbert-Schmidt operators from $\mU_i$ to $\mH$ with the Hilbert-Schmidt norm $\|\cdot\|_{\cL_2(\mU_i,\mH)}, i=1,2$.

\subsection{Multiscale stochastic Burgers equations with reflection}\label{msbewr}

In this subsection, we introduce the definition of the solutions for multiscale stochastic Burgers equations with reflection.

Let $A$ be the Laplace operator on $\mathbb{H}$
$$
A x:=\frac{\partial^2  x(\xi)}{\partial \xi^2}, \quad x \in D(A)=\mH^2([0,1]) \cap \mV,
$$
where $\mH^2([0,1])$ stands for the Sobolev space of functions defined on $[0,1]$ whose derivatives up to order $2$ belong to $\mH$. Let $\left\{e_k(\xi):=\sqrt{2} \sin (k \pi \xi)\right\}_{k \geq 1}$ be an orthonormal basis of $\mathbb{H}$ consisting of the eigenvectors of $A$, i.e.
$$
A e_k=-\lambda_k e_k \quad \text { with } \lambda_k=k^2 \pi^2.
$$

Define the bilinear operator $B(x, y): \mathbb{V} \times \mathbb{V} \rightarrow \mathbb{V}^*$ by
$$
B(x, y):=x \cdot \partial_{\xi} y,
$$
and the trilinear operator $b(x,y,z): \mV\times\mV\times\mV \rightarrow\mR$ by
$$
b(x,y,z):=\int_0^1x(\xi)\p_\xi y(\xi)z(\xi)\dif \xi={_{\mV^*}}\<B(x, y),z\>_\mV.
$$
We collect the properties about $B$ and $b$ in the Appendix.

With the above notations, the system (\ref{0eq}) can be rewritten as 
\be\left\{\begin{array}{ll}
\dif X_t^{\e,\g}=[AX_t^{\e,\g}+B(X_t^{\e,\g},X_t^{\e,\g})+F_1(X_t^{\e,\g},Y_t^{\e,\g})]\dif t+\sqrt\e G_1(X_t^{\e,\g})\dif W^1_t+\dif K^{1,\e,\g}_t,\\
\dif Y_t^{\e,\g}=\frac{1}{\g}[AY_t^{\e,\g}+F_2(X_t^{\e,\g},Y_t^{\e,\g})]\dif t+\frac{1}{\sqrt\g} G_2(X_t^{\e,\g},Y_t^{\e,\g})\dif W^2_t+\dif K^{2,\e,\g}_t,\\
X_t^{\e,\g}(\xi)\geq 0, \quad Y_t^{\e,\g}(\xi)\geq 0, \quad \xi\in[0,1],\\
X_t^{\e,\g}(0)=X_t^{\e,\g}(1)=Y_t^{\e,\g}(0)=Y_t^{\e,\g}(1)=0,\\
X_0^{\e,\g}=x_0, \quad Y_0^{\e,\g}=y_0,
\end{array}
\right.
\label{1eq}
\ee
where 
$$
A: \mV\to \mV^*, \quad B: \mV \times \mathbb{V} \rightarrow \mathbb{V}^*, \quad F_1: \mH\times\mH\to\mH, \quad G_1: \mH\to\cL_2(\mU_1,\mH),
$$
and
$$
F_2: \mH\times\mH\to\mH, \quad G_2: \mH\times\mH\to\cL_2(\mU_2,\mH),
$$
are measurable mappings, $W^1$ and $W^2$ are mutually independent $\mU_1$ and $\mU_2$-valued cylindrical Wiener processes on a complete filtered probability space $(\Omega, \sF, \{\sF_t\}_{t\geq 0}, \mP)$ respectively.

In the following, we present the definition of solutions for the system (\ref{1eq}) (cf. \cite{zhangt}).

\bd\label{soludefi}
$(X^{\e,\g}, K^{1,\e,\g}, Y^{\e,\g}, K^{2,\e,\g})$ is said to be a solution of the system (\ref{1eq}) if

$(i)$ $X_t^{\e,\g}, Y_t^{\e,\g}$ are $\mV$-valued $\sF_t$-measurable for any $t\geq 0$ and $X_t^{\e,\g}(\xi)\geq 0, Y_t^{\e,\g}(\xi)\geq 0$ a.e. for any $(t,\xi)\in\mR_+\times[0,1]$;

$(ii)$ $K^{1,\e,\g}, K^{2,\e,\g}$ are two random measures on $\mR_+\times[0,1]$ such that 

$(a)$ $\mE\left[({\rm Var}(K^{1,\e,\g})([0,T]\times [0,1]))^2\right]<+\infty$, $\mE\left[({\rm Var}(K^{2,\e,\g})([0,T]\times[0,1]))^2\right]<+\infty, \forall T \geq 0$, where ${\rm Var}(K^{1,\e,\g})([0,T]\times[0,1])$ denotes the total variation of $K^{1,\e,\g}$ on $[0, T]\times[0,1]$ defined by 
\ce
{\rm Var}(K^{1,\e,\g})([0,T]\times[0,1]):=\sup\limits_{\pi}\sum_{i=1}^n |K^{1,\e,\g}(E_i)|,
\de
and the supremum is taken over all partitions $\pi$ of the domain $[0,T]\times[0,1]$,

$(b)$ $K^{1,\e,\g}, K^{2,\e,\g}$ are adapted in the sense that for any bounded measurable mapping $\varphi$:
\ce
\int_0^t\int_0^1\varphi(s,\xi)K^{1,\e,\g}(\dif s, \dif \xi), \int_0^t\int_0^1\varphi(s,\xi)K^{2,\e,\g}(\dif s, \dif \xi) ~\mbox{are}~\sF_t\mbox{-measurable};
\de

$(iii)$ $(X^{\e,\g}, K^{1,\e,\g}, Y^{\e,\g}, K^{2,\e,\g})$ satisfies the parabolic stochastic partial differential equations with reflection in the following sense: $\forall t \in \mathbb{R}_{+}, \psi \in C_c^2([0,1])$ with $\psi(0)=\psi(1)=0$,
\ce
&&\<X_t^{\e,\g},\psi\>=\<x_0,\psi\>+\int_0^t\[{_{\mV^*}}\<AX_s^{\e,\g}+B(X_s^{\e,\g},X_s^{\e,\g}),\psi\>_\mV+\<F_1(X_s^{\e,\g},Y_s^{\e,\g}),\psi\>_\mH\]\dif s\\
&&\qquad\qquad\qquad+\sqrt\e \int_0^t \<\psi, G_1(X_s^{\e,\g})\dif W^1_s\>_\mH+\int_0^t\int_0^1\psi(\xi) \dif K^{1,\e,\g}(\dif s,\dif \xi), ~a.s.,\\
&&\<Y_t^{\e,\g},\psi\>=\<y_0,\psi\>+\frac{1}{\g}\int_0^t\[{_{\mV^*}}\<AY_s^{\e,\g},\psi\>_\mV+\<F_2(X_s^{\e,\g},Y_s^{\e,\g}),\psi\>_\mH\]\dif s\\
&&\qquad\qquad\qquad+\frac{1}{\sqrt\g}\int_0^t \<\psi,G_2(X_s^{\e,\g},Y_s^{\e,\g})\dif W^2_s\>_\mH+\int_0^t\int_0^1\psi(\xi) \dif K^{2,\e,\g}(\dif s,\dif \xi), ~a.s.;
\de

$(iv)$ for any $T>0$, $\int_0^T\int_0^1X_t^{\e,\g}(\xi)\dif K^{1,\e,\g}(\dif t, \dif \xi)=0, \quad \int_0^T\int_0^1Y_t^{\e,\g}(\xi)\dif K^{2,\e,\g}(\dif t, \dif \xi)=0$.
\ed

\subsection{LDPs}

In this subsection, we introduce the LDP.

Let $(\mX, \rho_\mX)$ be a Polish space. Let $\{X^{\e}, \e>0\}$ be a family of $\mX$-valued random variables defined on $(\Omega, \mathscr{F}, \{\mathscr{F}_t\}_{t\geq 0}, \mP)$.

\bd\label{rfde} 
$(i)$ A function $\Lambda: \mX\rightarrow[0,+\infty]$ is called a rate function on $\mX$, if for all $M\geq 0$, $\{\phi\in\mX: \Lambda(\phi)\leq M\}$ is a closed subset of $\mX$.

$(ii)$ A function $\Lambda: \mX\rightarrow[0,+\infty]$ is called a good rate function on $\mX$, if for all $M\geq 0$, $\{\phi\in\mX: \Lambda(\phi)\leq M\}$ is a compact subset of $\mX$.
\ed

\bd
We say that $\{X^{\e}\}$ satisfies the LDP with the speed $\e$ and the good rate function $\Lambda$, if for any subset $D\in \sB(\mX)$,
$$
-\inf\limits_{\phi\in {\rm Int}(D)}\Lambda(\phi)\leq\liminf_{\e\rightarrow 0}\e\log\mP(X^{\e}\in{\rm Int}(D))\leq \limsup\limits_{\e\rightarrow 0}\e\log\mP(X^{\e}\in \bar{D})\leq -\inf\limits_{\phi\in \bar{D}}\Lambda(\phi),
$$
where ${\rm Int}(D)$ and $\bar{D}$ denote the interior and the closure of $D$, respectively and they are taken in $\mX$.
\ed

\bd
We say that $\{X^{\e}\}$ satisfies the Laplace principle with the speed $\e$ and the good rate function $\Lambda$,  if for any real bounded continuous function $\Psi$ on $\mathbb{X}$,
\ce
\lim\limits_{\varepsilon\rightarrow 0}\varepsilon \log \mE\left\{\exp\left[-\frac{\Psi(X^{\e})}{\e}\right]\right\}=-\inf\limits_{\phi\in \mathbb{X}}\(\Psi(\phi)+\Lambda(\phi)\).
\de
\ed

In the sequel, we fix $T>0$ and take $\mX=C([0,T],\mH)$. Then $\mX$ is a Polish space. If $\Lambda$ is a good rate function, the LDP is equivalent to the Laplace principle (\cite{de}). Therefore, in order to obtain the LDP for $\{X^{\e}\}$, we state the conditions 
under which the Laplace principle holds. Set for any $N>0$, 
$$
\mathbf{D}^N=\left\{h\in L^2([0,T], \sU): \int_{0}^{T}|h(t)|_\sU^{2}\dif t \leq N\right\},
$$
where $\sU=\mU_1\times\mU_2$ is the product Hilbert space, and we equip $\mathbf{D}^{N}$ with the weak convergence topology in $L^2([0, T], \sU)$. So, $\mathbf{D}^{N}$ is metrizable as a compact Polish space. Let $\cA$ be  the collection of $(\sF_t)_{t\in[0,T]}$-predictable square integrable $\sU$-valued processes and $\mathcal{A}^N$ be the space of $\mathbf{D}^N$-valued random controls:
$$
\cA^N=\left\{u\in \mathcal{A}: u(\cdot, \omega) \in \mathbf{D}^N, \mathbb{P}\text { a.s. } \omega\right\}.
$$

Besides, let $\{e'_k, k\in\mN\}$ be the orthonormal basis on $\sU$. Define an $\sU$-valued cylindrical Wiener process $W$ as follows:
$$
W_t:=\sum_{k=1}^\infty \b_k(t)e'_k,
$$
where $\{\b_k, k\in\mN\}$ is a sequence of independent $1$-dimensional Brownian motions. So, there exists a Hilbert space $\sU_0$ such that $\sU\subset\sU_0$ is Hilbert-Schmidt and the paths of $W$ take values in $C([0,T],\sU_0)$. For cylindrical Wiener processes $W^1, W^2$, we can choose the projection operators $\Pi_1: \sU\to\mU_1$, $\Pi_2: \sU\to\mU_2$ such that
\ce
W_t^1=\Pi_1W_t, \quad W_t^2=\Pi_2W_t.
\de

Let $\mathcal{G}^\e: C([0,T],\sU_0)\rightarrow \mX$ be a measurable mapping for any $0<\e<1$.

\bco\label{cond}
There exists a measurable mapping $\mathcal{G}^0: C([0,T],\sU_0)\rightarrow \mX$ such that two following conditions hold:

$(i)$ For each $N<\infty$, any $\{h_\e\}\subset \mathbf{D}^N$ and $h\in \mathbf{D}^N$ satisfying $h_\e\to h$ as $\e\rightarrow0$, 
$$
\lim\limits_{\e\to0}\sup\limits_{t\in[0,T]}\left|\cG^0\left(\int_{0}^{\cdot}h_\e(s)\dif s\right)(t)-\cG^{0}\left(\int_{0}^{\cdot}h(s)\dif s\right)(t)\right|_\mH=0.
$$

$(ii)$ For each $N<\infty$, any $\{u_\e\}\subset \cA^N$ and any $\eta>0$,
$$
\lim _{\varepsilon \rightarrow 0}\mathbb{P}\left(\sup\limits_{t\in[0,T]}\left|\cG^{\varepsilon}\left(W+\frac{1}{\sqrt{\varepsilon} }\int_0^{\cdot} u_\e(s) \dif s\right)(t)-\cG^0\left(\int_0^{\cdot} u_\e(s)\dif s\right)(t)\right|_\mH>\eta\right)=0.
$$
\eco

The following result is due to \cite[Theorem 3.2]{aw}.

\bt\label{ldpth}
Set $X^\e=\mathcal{G}^\e\left(W\right)$. Suppose that Condition \ref{cond} holds. Then, the family $\{X^{\e}, \e>0\}$ satisfies the LDP on $\mX$, with the good rate function $\Lambda$ given by
\ce
\Lambda(\phi)=\inf\limits _{\substack{\{h\in L^2([0,T],\sU): \\
\phi=\cG^0(\int_0^{\cdot}h(s)\dif s)\}}}\left\{\frac{1}{2}\int_0^T|h(t)|_{\sU}^2\dif t\right\}, 
\de
with the convention $\inf\{\emptyset\}=\infty$.
\et

\section{Main result}\label{main}

In this section, we formulate the main result in this paper.

We recall system (\ref{1eq}), i.e.
\ce\left\{\begin{array}{ll}
\dif X_t^{\e,\g}=[AX_t^{\e,\g}+B(X_t^{\e,\g},X_t^{\e,\g})+F_1(X_t^{\e,\g},Y_t^{\e,\g})]\dif t+\sqrt\e G_1(X_t^{\e,\g})\dif W^1_t+\dif K^{1,\e,\g}_t,\\
\dif Y_t^{\e,\g}=\frac{1}{\g}[AY_t^{\e,\g}+F_2(X_t^{\e,\g},Y_t^{\e,\g})]\dif t+\frac{1}{\sqrt\g} G_2(X_t^{\e,\g},Y_t^{\e,\g})\dif W^2_t+\dif K^{2,\e,\g}_t,\\
X_t^{\e,\g}(\xi)\geq 0, \quad Y_t^{\e,\g}(\xi)\geq 0, \quad \xi\in[0,1],\\
X_t^{\e,\g}(0)=X_t^{\e,\g}(1)=Y_t^{\e,\g}(0)=Y_t^{\e,\g}(1)=0,\\
X_0^{\e,\g}=x_0, \quad Y_0^{\e,\g}=y_0.
\end{array}
\right.
\de

Assume:
\begin{enumerate}[$({\bf H}_{F_1, G_1})$]
\item There exist two constants $L_1, L_2\geq 0$ such that for any $x_1, x_2, y_1, y_2\in\mH$,
\ce
&&|F_1(x_1,y_1)-F_1(x_2,y_2)|_\mH\leq L_1(|x_1-x_2|_\mH+|y_1-y_2|_\mH),\\
&&\|G_1(x_1)-G_1(x_2)\|_{\cL_2(\mU_1,\mH)}\leq L_2|x_1-x_2|_\mH.
\de
\end{enumerate}
\begin{enumerate}[$({\bf H}^1_{F_2, G_2})$]
\item There exist four constants $L_3, L_4, L_5, L_6\geq 0$ such that for any $x_1, x_2, y_1, y_2\in\mH$,
\ce
&&|F_2(x_1,y_1)-F_2(x_2,y_2)|_\mH\leq L_3|x_1-x_2|_\mH+L_4|y_1-y_2|_\mH,\\
&&\|G_2(x_1,y_1)-G_2(x_2,y_2)\|_{\cL_2(\mU_2,\mH)}\leq L_5|x_1-x_2|_\mH+L_6|y_1-y_2|_\mH.
\de
\end{enumerate}
\begin{enumerate}[$({\bf H}^2_{G_2})$]
\item There exists a constant $L_7>0$ such that for any $x, y\in\mH$, $\varrho\in\mU_2$ and $\xi\in[0,1]$
\ce
|(G_2(x,y)\varrho)(\xi)|\geq L_7.
\de
\end{enumerate}
\begin{enumerate}[$({\bf H}^3_{F_2, G_2})$]
\item $L_4, L_6$ satisfy
\ce
\l_1-\frac{3}{2}L_4-L^2_6>0,
\de
where $\l_1>0$ is the smallest eigenvalue of $-A$.
\end{enumerate}
\begin{enumerate}[$({\bf H}^4_{G_2})$]
\item There exists a constant $L_8\geq L_7$ such that for any $x\in\mH$
\ce
\sup\limits_{y\in\mH}\|G_2(x,y)\|_{\cL_2(\mU_2,\mH)}\leq L_8(1+|x|_\mH).
\de
\end{enumerate}

\br
$(i)$ $({\bf H}_{F_1, G_1})$ implies that for $x, y\in\mH$
\be
|F_1(x,y)|_\mH\leq L'_1(1+|x|_\mH+|y|_\mH), \quad \|G_1(x)\|_{\cL_2(\mU_1,\mH)}\leq L'_2(1+|x|_\mH),
\label{f1g1linegrow}
\ee
where $L'_1, L'_2>0$ are two constants.

$(ii)$ $({\bf H}^1_{F_2, G_2})$ yields that for $x, y\in\mH$
\be
|F_2(x,y)|_\mH\leq C(1+|x|_\mH+|y|_\mH).
\label{f2g2linegrow}
\ee

$(iii)$ $({\bf H}^1_{F_2, G_2})$, $({\bf H}^3_{F_2, G_2})$, the Poincar\'e inequality and the Young inequality imply that for $x\in\mH, y\in\mV$
\be
&&2{_{\mV^*}}\<Ay, y\>_\mV+2\<F_2(x, y), y\>_\mH+\|G_2(x,y)\|^2_{\cL_2(\mU_2,\mH)}\no\\
&\leq&-2\|y\|_\mV^2+2\<F_2(x, y)-F_2(x, 0), y\>_\mH+2\<F_2(x, 0), y\>_\mH+2\|G_2(x,y)-G_2(x,0)\|^2_{\cL_2(\mU_2,\mH)}\no\\
&&+2\|G_2(x,0)\|^2_{\cL_2(\mU_2,\mH)}\no\\
&\leq&-2\l_1|y|_\mH^2+2L_4|y|_\mH^2+L_4|y|_\mH^2+2L^2_6|y|_\mH^2+C(1+|x|_\mH^2)\no\\
&=&-2\a |y|_\mH^2+C(1+|x|_\mH^2),
\label{invameasex}
\ee
where $\a:=\l_1-\frac{3}{2}L_4-L^2_6>0$.
\er

\br
$({\bf H}^1_{F_2, G_2})$ and $({\bf H}^2_{G_2})$ yield the existence and uniqueness of the invariant measure for the frozen equation, while $({\bf H}_{F_1, G_1})$, $({\bf H}^1_{F_2, G_2})$, $({\bf H}^2_{G_2})$, $({\bf H}^3_{F_2, G_2})$ and $({\bf H}^4_{G_2})$ imply the LDP for the system (\ref{1eq}).
\er

Based on Theorem 2.3 in \cite{lrsx}, following the idea of Theorem 2.1 in \cite{my}, we can obtain the following well-posedness result, whose proof is omitted.

\bt\label{wellpose}
Assume that $({\bf H}_{F_1, G_1})$, $({\bf H}^1_{F_2, G_2})$, $({\bf H}^2_{G_2})$ and $({\bf H}^3_{F_2, G_2})$ hold. Then for any $x_0, y_0\in\mV, x_0(\xi)\geq 0, y_0(\xi)\geq 0$, the system (\ref{1eq}) has a unique solution $(X^{\e,\g}, K^{1,\e,\g},Y^{\e,\g}, \\ K^{2,\e,\g})$ such that 
$$
\mE\left[\sup\limits_{t\in[0,T]}|X^{\e,\g}_t|_\mH^2+\int_0^T\|X^{\e,\g}_t\|^2_\mV\dif t\right]<\infty, \quad \mE\left[\sup\limits_{t\in[0,T]}|Y^{\e,\g}_t|_\mH^2+\int_0^T\|Y^{\e,\g}_t\|^2_\mV\dif t\right]<\infty.
$$
\et

Next, fix $x\in\mH$ and consider the following SPDE with reflection:
\be\left\{\begin{array}{ll}
\dif Y_t^{x,y_0}=[AY_t^{x,y_0}+F_2(x,Y_t^{x,y_0})]\dif t+G_2(x,Y_t^{x,y_0})\dif \tilde W^2_t+\dif K^{2,x,y_0}_t,\\
Y_t^{x,y_0}(\xi)\geq 0,\quad \xi\in[0,1],\\
Y_t^{x,y_0}(0)=Y_t^{x,y_0}(1)=0,\\
Y_0^{x,y_0}=y_0,
\end{array}
\right.
\label{frozequa}
\ee
where $\tilde W^2$ is an $\mU_2$-valued cylindrical Wiener process independent of $W^1, W^2$. Under $({\bf H}^1_{F_2,G_2})$, Eq.(\ref{frozequa}) has a unique solution $(Y^{x,y_0}, K^{2,x,y_0})$ such that $\mE[\sup\limits_{t\in[0,T]}|Y^{x,y_0}_t|_\mH^2+\int_0^T\|Y^{x,y_0}_t\|^2_\mV\dif t]<\infty$ (\cite{xz}). Moreover, we notice that $Y^{x,y_0}$ is a Markov process and then has a unique invariant probability measure $\mu^x$ under $({\bf H}^2_{G_2})$ (\cite{yz}).

Set
\ce
\bar F_1(x):=\int_{\mH}F_1(x,y)\mu^x(\dif y),
\de
and we construct the following PDE with reflection: for any $h\in L^2([0,T],\sU)$
\be\left\{\begin{array}{ll}
\dif \bar X_t^{h}=[A\bar X_t^{h}+B(\bar X_t^{h},\bar X_t^{h})+\bar F_1(\bar X_t^{h})]\dif t+G_1(\bar X_t^{h})\Pi_1h(t)\dif t+\dif \bar K^{1,h}_t,\\
\bar X_t^{h}(\xi)\geq 0, \quad \xi\in[0,1],\\
\bar X_t^{h}(0)=\bar X_t^{h}(1)=0,\\
\bar X_0^{h}=x_0.
\end{array}
\right.
\label{e0conteq}
\ee
By Lemma \ref{exisuniq}, Eq.(\ref{e0conteq}) has a unique solution $(\bar{X}^{h}, \bar{K}^{1,h})$ such that $\sup\limits_{t\in[0,T]}|\bar{X}^{h}_t|_\mH^2+\int_0^T\|\bar{X}^{h}_t\|^2_\mV\dif t<\infty$. 

Now we state the main result in this paper.

\bt\label{fldpth}
Suppose that $({\bf H}_{F_1, G_1})$, $({\bf H}^1_{F_2, G_2})$, $({\bf H}^2_{G_2})$, $({\bf H}^3_{F_2, G_2})$ and $({\bf H}^4_{G_2})$ hold. If 
\ce
\lim\limits_{\e\rightarrow0}\frac{\g}{\e}=0,
\de
then $\{X^{\e,\g}, \e>0\}$ satisfies the LDP in $C([0,T], \mH)$ with the rate function $\Lambda$ given by
\ce
\Lambda(\phi)=\inf\limits _{\substack{\{h\in L^2([0,T],\sU): \\\phi=\bar{X}^{h}\}}}\left\{\frac{1}{2}\int_0^T|h(t)|_{\sU}^2\dif t\right\}, \quad \phi\in C([0,T], \mH),
\de
with the convention $\inf\{\emptyset\}=\infty$.
\et

The proof of Theorem \ref{fldpth} is placed in the next section.

\section{Proof of Theorem \ref{fldpth}}\label{fldpthproo}

In this section, we prove Theorem \ref{fldpth} through Theorem \ref{ldpth}. First of all, we estimate the solution of the frozen equation (\ref{frozequa}) in Subsection \ref{frozequaesti}. Next, in Subsection \ref{vericond1} Condition \ref{cond} $(i)$ is verified. Then we justify Condition \ref{cond} $(ii)$ in Subsection \ref{vericond2}. Finally, we finish the proof of Theorem \ref{fldpth}.

\subsection{Some estimates about the frozen equation (\ref{frozequa})}\label{frozequaesti}

In this subsection, we present some results about the frozen equation (\ref{frozequa}). We begin with some key estimates.

\bl
Under $({\bf H}^1_{F_2,G_2})$ and $({\bf H}^3_{F_2,G_2})$ it holds that for any $x, x_1, x_2, y_1, y_2\in\mH$ and $t\geq 0$
\be
&&\mE|Y_t^{x,y_0}|_\mH^2\leq C(1+|x|_\mH^2)+|y_0|_\mH^2 e^{-2\a t},\label{frozmomeesti}\\
&&\mE|Y_t^{x_1,y_1}-Y_t^{x_2,y_2}|_\mH^2\leq C|x_1-x_2|_\mH^2+|y_1-y_2|_\mH^2 e^{-2\a t}. \label{frozdiffesti}
\ee
\el
\begin{proof}
First of all, by the It\^o formula, it holds that for any $l_1>0$ 
\ce
e^{l_1t}|Y_t^{x,y_0}|_\mH^2&=&|y_0|_\mH^2+l_1\int_0^te^{l_1s}|Y_s^{x,y_0}|_\mH^2\dif s+2\int_0^{t} e^{l_1s}{_{\mV^*}}\<AY_s^{x,y_0},Y_s^{x,y_0}\>_{\mV}\dif s\\
&&+2\int_0^te^{l_1s}\<F_2(x,Y_s^{x,y_0}),Y_s^{x,y_0}\>_\mH\dif s+2\int_0^te^{l_1s}\<Y_s^{x,y_0}, G_2(x,Y_s^{x,y_0})\dif \tilde W^2_s\>_\mH\\
&&+2\int_0^t\int_0^1e^{l_1s}Y_s^{x,y_0}(\xi)K^{2,x,y_0}(\dif s,\dif \xi)+\int_0^te^{l_1s}\|G_2(x,Y_s^{x,y_0})\|^2_{\cL_2(\mU_2,\mH)}\dif s.
\de
Then (\ref{invameasex}) implies that
\ce
&&2{_{\mV^*}}\<AY_s^{x,y_0},Y_s^{x,y_0}\>_{\mV}+2\<F_2(x,Y_s^{x,y_0}),Y_s^{x,y_0}\>_\mH+\|G_2(x,Y_s^{x,y_0})\|^2_{\cL_2(\mU_2,\mH)}\\
&\leq& -2\a |Y_s^{x,y_0}|_\mH^2+C(1+|x|_\mH^2).
\de
From Definition \ref{soludefi}, it follows that 
$$
2\int_0^t\int_0^1e^{l_1 s}Y_s^{x,y_0}(\xi)K^{2,x,y_0}(\dif s,\dif \xi)=0.
$$
Collecting the above deduction and taking the expectation on two sides of the above equality, we have that
\ce
e^{l_1t}\mE|Y_t^{x,y_0}|_\mH^2\leq |y_0|_\mH^2+l_1\int_0^te^{l_1s}\mE|Y_s^{x,y_0}|_\mH^2\dif s-2\a\int_0^{t} e^{l_1s}\mE|Y_s^{x,y_0}|_\mH^2\dif s+C(1+|x|_\mH^2)\int_0^te^{l_1s}\dif s.
\de
 So, we take $l_1=2\a$ and obtain that
 \ce
\mE|Y_t^{x,y_0}|_\mH^2\leq C(1+|x|_\mH^2)+|y_0|_\mH^2e^{-2\a t}.
 \de

Next, we assume that $(Y^{x_1,y_1}, K^{2,x_1,y_1})$ and $(Y^{x_2,y_2}, K^{2,x_2,y_2})$ are the solutions of Eq.(\ref{frozequa}) with replacing $(x, y_0)$ by $(x_1, y_1), (x_2, y_2)$, respectively. The It\^o formula implies that for any $l_2>0$
\ce
e^{l_2 t}|Y_t^{x_1,y_1}-Y_t^{x_2,y_2}|_\mH^2&=&|y_1-y_2|_\mH^2+l_2\int_0^te^{l_2 s}|Y_s^{x_1,y_1}-Y_s^{x_2,y_2}|_\mH^2\dif s\\
&&+2\int_0^te^{l_2 s}{_{\mV^*}}\<A(Y_s^{x_1,y_1}-Y_s^{x_2,y_2}),Y_s^{x_1,y_1}-Y_s^{x_2,y_2}\>_{\mV}\dif s\\
&&+2\int_0^te^{l_2 s}\<F_2(x_1,Y_s^{x_1,y_1})-F_2(x_2,Y_s^{x_2,y_2}),Y_s^{x_1,y_1}-Y_s^{x_2,y_2}\>_\mH\dif s\\
&&+2\int_0^te^{l_2 s}\<Y_s^{x_1,y_1}-Y_s^{x_2,y_2}, \(G_2(x_1,Y_s^{x_1,y_1})-G_2(x_2,Y_s^{x_2,y_2})\)\dif \tilde W^2_s\>_\mH\\
&&+2\int_0^t\int_0^1e^{l_2 s}(Y_s^{x_1,y_1}(\xi)-Y_s^{x_2,y_2}(\xi))K^{2,x_1,y_1}(\dif s,\dif \xi)\\
&&-2\int_0^t\int_0^1e^{l_2 s}(Y_s^{x_1,y_1}(\xi)-Y_s^{x_2,y_2}(\xi))K^{2,x_2,y_2}(\dif s,\dif \xi)\\
&&+\int_0^te^{l_2 s}\|G_2(x_1,Y_s^{x_1,y_1})-G_2(x_2,Y_s^{x_2,y_2})\|^2_{\cL_2(\mU_2,\mH)}\dif s.
\de
By the Poincar\'e inequality, the H\"older inequality and $({\bf H}^1_{F_2,G_2})$, it holds that
\ce
&&2{_{\mV^*}}\<A(Y_s^{x_1,y_1}-Y_s^{x_2,y_2}),Y_s^{x_1,y_1}-Y_s^{x_2,y_2}\>_{\mV}\leq -2\l_1|Y_s^{x_1,y_1}-Y_s^{x_2,y_2}|_\mH^2,\\
&&2\<F_2(x_1,Y_s^{x_1,y_1})-F_2(x_2,Y_s^{x_2,y_2}),Y_s^{x_1,y_1}-Y_s^{x_2,y_2}\>_\mH\\
&\leq&2\<F_2(x_1,Y_s^{x_1,y_1})-F_2(x_1,Y_s^{x_2,y_2}),Y_s^{x_1,y_1}-Y_s^{x_2,y_2}\>_\mH\\
&&+2\<F_2(x_1,Y_s^{x_2,y_2})-F_2(x_2,Y_s^{x_2,y_2}),Y_s^{x_1,y_1}-Y_s^{x_2,y_2}\>_\mH\\
&\leq& 2L_4|Y_s^{x_1,y_1}-Y_s^{x_2,y_2}|_\mH^2+L_4|Y_s^{x,y_1}-Y_s^{x,y_2}|_\mH^2+C|x_1-x_2|_\mH^2,
\de
and
\ce
&&\|G_2(x_1,Y_s^{x_1,y_1})-G_2(x_2,Y_s^{x_2,y_2})\|^2_{\cL_2(\mU_2,\mH)}\\
&\leq&2\|G_2(x_1,Y_s^{x_1,y_1})-G_2(x_1,Y_s^{x_2,y_2})\|^2_{\cL_2(\mU_2,\mH)}\\
&&+2\|G_2(x_1,Y_s^{x_2,y_2})-G_2(x_2,Y_s^{x_2,y_2})\|^2_{\cL_2(\mU_2,\mH)}\\
&\leq&2L_6^2|Y_s^{x_1,y_1}-Y_s^{x_2,y_2}|_\mH^2+2L_5^2|x_1-x_2|_\mH^2.
\de
Note that $\int_0^t\int_0^1 e^{l_2 s}Y_s^{x_1,y_1}(\xi)K^{2,x_1,y_1}(\dif s,\dif \xi)=0$ and $Y_s^{x_2,y_2}(\xi)\geq 0$. Thus, we infer that
\ce
&&2\int_0^t\int_0^1e^{l_2 s}\(Y_s^{x_1,y_1}(\xi)-Y_s^{x_2,y_2}(\xi)\)K^{2,x_1,y_1}(\dif s,\dif \xi)\\
&=&2\int_0^t\int_0^1e^{l_2 s}Y_s^{x_1,y_1}(\xi)K^{2,x_1,y_1}(\dif s,\dif \xi)-2\int_0^t\int_0^1e^{l_2 s}Y_s^{x_2,y_2}(\xi)K^{2,x_1,y_1}(\dif s,\dif \xi)\\
&\leq& 0.
\de
From the same deduction to that for the above inequality, it follows that
\ce
-2\int_0^t\int_0^1e^{l_2 s}(Y_s^{x_1,y_1}(\xi)-Y_s^{x_2,y_2}(\xi))K^{2,x_2,y_2}(\dif s,\dif \xi)\leq 0.
\de
Combining the above deduction and taking the expectation on two sides, we conclude that
\ce
e^{l_2 t}\mE|Y_t^{x_1,y_1}-Y_t^{x_2,y_2}|_\mH^2&\leq& |y_1-y_2|_\mH^2+l_2\int_0^te^{l_2 s}\mE|Y_s^{x_1,y_1}-Y_s^{x_2,y_2}|_\mH^2\dif s\\
&&-2\a\int_0^te^{l_2 s}\mE|Y_s^{x_1,y_1}-Y_s^{x_2,y_2}|_\mH^2\dif s+C|x_1-x_2|_\mH^2\int_0^t e^{l_2 s}\dif s.
\de
Set $l_2=2\a$ and it holds that
\ce
\mE|Y_t^{x_1,y_1}-Y_t^{x_2,y_2}|_\mH^2\leq C|x_1-x_2|_\mH^2+|y_1-y_2|_\mH^2 e^{-2\a t}.
\de
The proof is complete.
\end{proof}

By (\ref{frozmomeesti}) and the definition of $\mu^{x}$, it holds that for any $t\geq 0$
\ce
\int_{\mH}|y|_\mH^{2}\mu^{x}(\dif y)
&=&\int_{\mH}\mE|Y_{t}^{x,y}|_\mH^{2}\mu^{x}(\dif y)\leq \int_{\mH}\(|y|_\mH^{2}e^{-2\a t}+C(1+|x|_\mH^{2})\)\mu^{x}(\dif y)\no\\
&=& e^{-2\a t}\int_{\mH}|y|_\mH^{2}\mu^{x}(\dif y)+C(1+|x|_\mH^{2}),
\de
and furthermore
\be
\int_{\mH}|y|_\mH^{2}\mu^{x}(\dif y)\leq C(1+|x|_\mH^{2}).
\label{inu2}
\ee

\bl\label{emb1}
 Suppose that $({\bf H}_{F_1, G_1})$, $({\bf H}^1_{F_2,G_2})$, $({\bf H}^2_{G_2})$ and $({\bf H}^3_{F_2, G_2})$ hold. Then, for any $t\geq 0$ there exists a constant $C>0$ such that for $x,y\in\mH$
\be
|\mE F_{1}(x,Y_{t}^{x,y})-\bar{F}_{1}(x)|_\mH^{2}\leq Ce^{-2\a t}(1+|x|_\mH^{2}+|y|_\mH^{2}).
\label{meu2}
\ee
\el
\begin{proof}
Based on simple calculations, one can obtain that
\ce
&&|\mE F_{1}(x,Y_{t}^{x,y})-\bar{F}_{1}(x)|_\mH^{2}
=\Big|\mE F_{1}(x,Y_{t}^{x,y})-\int_{\mH}F_{1}(x,z)\mu^{x}(\dif z)\Big|_\mH^{2}\\
&=&\Big|\mE F_{1}(x,Y_{t}^{x,y})-\int_{\mH}\mE F_{1}(x,Y_{t}^{x,z})\mu^{x}(\dif z)\Big|_\mH^{2}\\
&\leq&\int_{\mH}\mE |F_{1}(x,Y_{t}^{x,y})- F_{1}(x,Y_{t}^{x,z})|_\mH^{2}\mu^{x}(\dif z)\\
&\leq& L^2_1\int_{\mH}\mE |Y_{t}^{x,y}-Y_{t}^{x,z}|_\mH^{2}\mu^{x}(\dif z)\\
&\leq& L^2_1\int_{\mH}|y-z|_\mH^{2}e^{-2\a t}\mu^{x}(\dif z)\\
&\leq&2L^2_1e^{-2\a t}\left(|y|_\mH^{2}+\int_{\mH}|z|_\mH^{2}\mu^{x}(\dif z)\right),
\de
where the second and third inequalities are based on $({\bf H}_{F_1, G_1})$ and (\ref{frozdiffesti}), respectively. Finally, (\ref{inu2}) implies the required estimate.
\end{proof}

\subsection{Verification of Condition \ref{cond} $(i)$}\label{vericond1}

In this subsection, we justify Condition \ref{cond} $(i)$. 

Firstly, the well-posedness of Eq.(\ref{e0conteq}) is determined in the following lemma.

\bl\label{exisuniq}
Suppose that $({\bf H}_{F_1, G_1})$, $({\bf H}^1_{F_2,G_2})$, $({\bf H}^2_{G_2})$ and $({\bf H}^3_{F_2, G_2})$ hold. Then Eq.(\ref{e0conteq}) has a unique solution $(\bar{X}^{h}, \bar{K}^{1,h})$ such that $\bar{X}^{h}$ belongs to $C([0,T],\mH)\cap L^2([0,T],\mV)$.
\el
\begin{proof}
First of all, we prove that $\bar{F}_1$ is Lipschitz continuous in $x$. Indeed, from (\ref{meu2}) and (\ref{frozdiffesti}), it follows that for any $x_1, x_2\in\mH$ and $t\geq 0$
\ce
|\bar{F}_1(x_1)-\bar{F}_1(x_2)|_\mH^2&\leq& 3|\bar{F}_1(x_1)-\mE F_{1}(x_1,Y_{t}^{x_1,y_0})|_\mH^2\\
&&+3|\mE F_{1}(x_1,Y_{t}^{x_1,y_0})-\mE F_{1}(x_2,Y_{t}^{x_2,y_0})|_\mH^2\\
&&+3|\mE F_{1}(x_2,Y_{t}^{x_2,y_0})-\bar{F}_1(x_2)|_\mH^2\\
&\leq& 3Ce^{-2\a t}(2+|x_1|_\mH^{2}+|x_2|_\mH^{2}+2|y_0|_\mH^{2})+3C|x_1-x_2|_\mH^2.
\de
Letting $t\rightarrow \infty$, we obtain that
\be
|\bar{F}_1(x_1)-\bar{F}_1(x_2)|_\mH^2\leq 3C|x_1-x_2|_\mH^2.
\label{barf1lip}
\ee

Next, for any $n\in\mN$, we consider the following penalized Burgers type equation associated with Eq.(\ref{e0conteq}):
\be\left\{\begin{array}{ll}
\dif \bar X_t^{h,n}=[A\bar X_t^{h,n}+B(\bar X_t^{h,n},\bar X_t^{h,n})+\bar F_1(\bar X_t^{h,n})]\dif t+G_1(\bar X_t^{h,n})\Pi_1h(t)\dif t+n\bar X_t^{h,n-}\dif t,\\
\bar X_t^{h,n}(0)=\bar X_t^{h,n}(1)=0,\\
\bar X_0^{h,n}=x_0,
\end{array}
\right.
\label{e0conteqapp}
\ee
where $\bar X_t^{h,n-}(\xi):=-\min\{\bar X_t^{h,n}(\xi),0\}$. Under (\ref{barf1lip}) and $({\bf H}_{F_1, G_1})$, by Corollary 7.7 in \cite{bm}, we can obtain that Eq.(\ref{e0conteqapp}) has a unique solution $\bar X^{h,n}$ such that 
\ce
\sup\limits_{t\in[0,T]}|\bar{X}_t^{h,n}|^2_\mH+\int_0^T\|\bar{X}_t^{h,n}\|^2_\mV\dif t<\infty.
\de

Finally, by leveraging $\bar{X}^{h,n}$ and following an argument analogous to that of Proposition 2.1 in \cite{wzz}, we deduce that Eq.(\ref{e0conteq}) admits a unique solution $(\bar{X}^{h}, \bar{K}^{1,h})$, where $\bar{X}^{h} \in C([0,T], \mH) \cap L^2([0,T], \mV)$ and 
\be
\lim\limits_{n\rightarrow\infty}\left\{\sup\limits_{t\in[0,T]}|\bar X_t^{h,n}-\bar{X}_t^{h}|_\mH^2+\int_0^T\|\bar X_t^{h,n}-\bar{X}_t^{h}\|_\mV^2\dif t\right\}=0.
\label{barxhnbarxh}
\ee
The proof is complete.
\end{proof}

Next, about $\bar{X}^{h,n}$ we have the following estimates (cf. \cite{wzz})
\be
&&\sup\limits_{n}\sup\limits_{h\in{\bf D}^N}\(\sup\limits_{t\in[0,T]}|\bar{X}_t^{h,n}|^2_\mH+\int_0^T\|\bar{X}_t^{h,n}\|^2_\mV\dif t\)<\infty,\label{unifbarxhn}\\
&&\sup\limits_{n}\sup\limits_{h\in{\bf D}^N}\(n\|\bar{X}^{h,n-}\|_{L^1([0,T]\times[0,1])}\)<\infty,\label{unifnbarxhn}\\
&&\sup\limits_{n}\sup\limits_{h\in{\bf D}^N}\sup\limits_{t\in[0,T]}\left(\exp\left\{-4\int_0^t\|\bar{X}_s^{h,n}\|^2_\mV\dif s\right\}\|\bar{X}_t^{h,n}\|^2_\mV\right)<\infty, \label{unifexpbarx}\\
&&\lim\limits_{n\rightarrow\infty}\sup\limits_{h\in{\bf D}^N}\sup\limits_{t\in[0,T]}|\bar{X}_t^{h,n-}|_\mH=0,\label{unifbarxhn4}
\ee
where $L^1([0,T]\times[0,1])$ denotes the Sobolev space of absolutely integrable functions on $[0,T]\times[0,1]$.

\bl\label{unihdn}
Assume that $({\bf H}_{F_1, G_1})$, $({\bf H}^1_{F_2,G_2})$, $({\bf H}^2_{G_2})$ and $({\bf H}^3_{F_2, G_2})$ hold. Then for any $N>0$ it holds that
\be
\lim\limits_{n\rightarrow\infty}\sup\limits_{h\in{\bf D}^N}\left\{\sup\limits_{t\in[0,T]}|\bar X_t^{h,n}-\bar{X}_t^{h}|_\mH^2+\int_0^T\|\bar X_t^{h,n}-\bar{X}_t^{h}\|_\mV^2\dif t\right\}=0.
\label{hnhdiff}
\ee
\el
\begin{proof}
If we prove that
\be
\lim\limits_{n,m\rightarrow\infty}\sup\limits_{h\in{\bf D}^N}\left\{\sup\limits_{t\in[0,T]}|\bar X_t^{h,n}-\bar{X}_t^{h,m}|_\mH^2+\int_0^T\|\bar X_t^{h,n}-\bar{X}_t^{h,m}\|_\mV^2\dif t\right\}=0,
\label{hnhmdiff}
\ee
(\ref{barxhnbarxh}) implies (\ref{hnhdiff}). 

Next, we are devoted to showing (\ref{hnhmdiff}). First of all, set for $m, n\in\mN$
\ce
\Sigma_t:=\exp\left\{-\b_1\int_0^t(\|\bar X_s^{h,n}\|_\mV^2+\|\bar X_s^{h,m}\|_\mV^2)\dif s\right\}, \quad t\geq 0,
\de
where the constant $\b_1>0$ is determined later, and by the chain rule it holds that
\ce
\Sigma_t|\bar X_t^{h,n}-\bar{X}_t^{h,m}|_\mH^2&=&-\b_1\int_0^t(\|\bar X_s^{h,n}\|_\mV^2+\|\bar X_s^{h,m}\|_\mV^2)\Sigma_s|\bar X_s^{h,n}-\bar{X}_s^{h,m}|_\mH^2\dif s\\
&&+2\int_0^t\Sigma_s{_{\mV^*}}\<A(\bar X_s^{h,n}-\bar{X}_s^{h,m}),\bar X_s^{h,n}-\bar{X}_s^{h,m}\>_\mV\dif s\\
&&+2\int_0^t\Sigma_s{_{\mV^*}}\<B(\bar X_s^{h,n},\bar X_s^{h,n})-B(\bar X_s^{h,m},\bar X_s^{h,m}),\bar X_s^{h,n}-\bar{X}_s^{h,m}\>_\mV\dif s\\
&&+2\int_0^t\Sigma_s\<\bar F_1(\bar X_s^{h,n})-\bar F_1(\bar X_s^{h,m}),\bar X_s^{h,n}-\bar{X}_s^{h,m}\>_\mH\dif s\\
&&+2\int_0^t\Sigma_s\<G_1(\bar X_s^{h,n})\Pi_1h(s)-G_1(\bar X_s^{h,m})\Pi_1h(s),\bar X_s^{h,n}-\bar{X}_s^{h,m}\>_\mH\dif s\\
&&+2\int_0^t\Sigma_s\<n\bar X_s^{h,n-},\bar X_s^{h,n}-\bar{X}_s^{h,m}\>_\mH\dif s\\
&&-2\int_0^t\Sigma_s\<m\bar X_s^{h,m-},\bar X_s^{h,n}-\bar{X}_s^{h,m}\>_\mH\dif s\\
&=:&I_1+I_2+I_3+I_4+I_5+I_6+I_7.
\de
For $I_2$, by the definition of $A$, we infer that
\ce
I_2=-2\int_0^t\Sigma_s\|\bar X_s^{h,n}-\bar{X}_s^{h,m}\|_\mV^2\dif s.
\de
For $I_3$, Lemma \ref{Bbprop3} and the Young inequality imply that
\ce
I_3&\leq&2\int_0^t\Sigma_s 2|\bar X_s^{h,n}-\bar{X}_s^{h,m}|_\mH(\|\bar X_s^{h,n}\|_\mV+\|\bar X_s^{h,m}\|_\mV)\|\bar X_s^{h,n}-\bar{X}_s^{h,m}\|_\mV\dif s\\
&\leq&\int_0^t\Sigma_s\|\bar X_s^{h,n}-\bar{X}_s^{h,m}\|_\mV^2\dif s+8\int_0^t\Sigma_s|\bar X_s^{h,n}-\bar{X}_s^{h,m}|_\mH^2(\|\bar X_s^{h,n}\|^2_\mV+\|\bar X_s^{h,m}\|^2_\mV)\dif s.
\de
For $I_4$, (\ref{barf1lip}) yields that
\ce
I_4\leq 2C\int_0^t\Sigma_s|\bar X_s^{h,n}-\bar{X}_s^{h,m}|_\mH^2\dif s.
\de
For $I_5$, from $({\bf H}_{F_1, G_1})$ and the H\"older inequality it follows that
\ce
I_5&\leq&2\int_0^t\Sigma_s|G_1(\bar X_s^{h,n})\Pi_1h(s)-G_1(\bar X_s^{h,m})\Pi_1h(s)|_\mH|\bar X_s^{h,n}-\bar{X}_s^{h,m}|_\mH\dif s\\
&\leq&2\sup\limits_{s\in[0,t]}\Sigma^{1/2}_s|\bar X_s^{h,n}-\bar{X}_s^{h,m}|_\mH \int_0^t\Sigma^{1/2}_s\|G_1(\bar X_s^{h,n})-G_1(\bar X_s^{h,m})\|_{\cL_2(\mU_1,\mH)}|h(s)|_\sU\dif s\\
&\leq&\frac{1}{2}\sup\limits_{s\in[0,t]}\Sigma_s|\bar X_s^{h,n}-\bar{X}_s^{h,m}|_\mH^2+2NL_1^2\int_0^t\Sigma_s|\bar X_s^{h,n}-\bar{X}_s^{h,m}|_\mH^2\dif s. 
\de
For $I_6$, noticing $\<n\bar X_s^{h,n-},\bar X_s^{h,n+}\>=0$ and $-\<n\bar X_s^{h,n-},\bar X_s^{h,n-}\>\leq 0, -\<n\bar X_s^{h,n-},\bar X_s^{h,m+}\>\leq 0$, we obtain that
\ce
I_6&=&2\int_0^t\Sigma_s\<n\bar X_s^{h,n-},\bar X_s^{h,n}-\bar{X}_s^{h,m}\>_\mH\dif s\\
&=&2\int_0^t\Sigma_s\<n\bar X_s^{h,n-},\bar X_s^{h,n+}-\bar X_s^{h,n-}-\bar X_s^{h,m+}+\bar X_s^{h,m-}\>_\mH\dif s\\
&\leq&2\int_0^t\Sigma_s\<n\bar X_s^{h,n-},\bar X_s^{h,m-}\>_\mH\dif s\leq 2\int_0^t \Sigma_s\|\bar X_s^{h,m-}\|_{L^\infty([0,1])}n\int_0^1\bar X_s^{h,n-}(\xi)\dif \xi\dif s\\
&\leq&2\(\sup\limits_{s\in[0,T]}\Sigma_s\|\bar X_s^{h,m-}\|_{L^\infty([0,1])}\) n\|\bar X^{h,n-}\|_{L^1([0,T]\times[0,1])},
\de
where $L^\infty([0,1])$ denotes the Sobolev space of essentially bounded functions on $[0,1]$. And similarly,
\ce
I_7\leq 2\(\sup\limits_{s\in[0,T]}\Sigma_s\|\bar X_s^{h,n-}\|_{L^\infty([0,1])}\) m\|\bar X^{h,m-}\|_{L^1([0,T]\times[0,1])}.
\de

Collecting the above deduction and taking $\b_1\geq 8$, we conclude that
\ce
&&\sup\limits_{t\in[0,T]}\Sigma_t|\bar X_t^{h,n}-\bar{X}_t^{h,m}|_\mH^2+\int_0^T\Sigma_s\|\bar X_s^{h,n}-\bar{X}_s^{h,m}\|_\mV^2\dif s\\
&\leq&4(C+NL_1^2)\int_0^T\sup\limits_{r\in[0,s]}\Sigma_r|\bar X_r^{h,n}-\bar{X}_r^{h,m}|_\mH^2\dif s\\
&&+4\(\sup\limits_{s\in[0,T]}\Sigma_s\|\bar X_s^{h,m-}\|_{L^\infty([0,1])}\) n\|\bar X^{h,n-}\|_{L^1([0,T]\times[0,1])}\\
&&+4\(\sup\limits_{s\in[0,T]}\Sigma_s\|\bar X_s^{h,n-}\|_{L^\infty([0,1])}\) m\|\bar X^{h,m-}\|_{L^1([0,T]\times[0,1])},
\de
which together with the Gronwall inequality implies that
\ce
&&\sup\limits_{t\in[0,T]}\Sigma_t|\bar X_t^{h,n}-\bar{X}_t^{h,m}|_\mH^2+\int_0^T\Sigma_s\|\bar X_s^{h,n}-\bar{X}_s^{h,m}\|_\mV^2\dif s\\
&\leq&C\bigg(\(\sup\limits_{s\in[0,T]}\Sigma_s\|\bar X_s^{h,m-}\|_{L^\infty([0,1])}\) n\|\bar X^{h,n-}\|_{L^1([0,T]\times[0,1])}\\
&&+\(\sup\limits_{s\in[0,T]}\Sigma_s\|\bar X_s^{h,n-}\|_{L^\infty([0,1])}\) m\|\bar X^{h,m-}\|_{L^1([0,T]\times[0,1])}\bigg),
\de
and by (\ref{unifnbarxhn})
\ce
&&\sup\limits_{h\in{\bf D}^N}\left\{\sup\limits_{t\in[0,T]}\Sigma_t|\bar X_t^{h,n}-\bar{X}_t^{h,m}|_\mH^2+\int_0^T\Sigma_s\|\bar X_s^{h,n}-\bar{X}_s^{h,m}\|_\mV^2\dif s\right\}\\
&\leq&C\Bigg(\sup\limits_{h\in{\bf D}^N}\sup\limits_{s\in[0,T]}\left(\exp\left\{-2\int_0^s\|\bar X_r^{h,m}\|_\mV^2\dif r\right\}\|\bar X_s^{h,m-}\|_{L^\infty([0,1])}\right)\\
&&\quad+\sup\limits_{h\in{\bf D}^N}\sup\limits_{s\in[0,T]}\left(\exp\left\{-2\int_0^s\|\bar X_r^{h,n}\|_\mV^2\dif r\right\}\|\bar X_s^{h,n-}\|_{L^\infty([0,1])}\right)\Bigg).
\de

Finally, we observe $\sup\limits_{h\in{\bf D}^N}\sup\limits_{s\in[0,T]}\left(\exp\left\{-2\int_0^s\|\bar X_r^{h,m}\|_\mV^2\dif r\right\}\|\bar X_s^{h,m-}\|_{L^\infty([0,1])}\right)$. By the Sobolev imbedding theorem, for any $\t>0$, there exists a constant $C_\t>0$ such that
\ce
&&\sup\limits_{h\in{\bf D}^N}\sup\limits_{s\in[0,T]}\left(\exp\left\{-4\int_0^s\|\bar X_r^{h,m}\|_\mV^2\dif r\right\}\|\bar X_s^{h,m-}\|^2_{L^\infty([0,1])}\right)\\
&\leq&\t\sup\limits_{h\in{\bf D}^N}\sup\limits_{s\in[0,T]}\left(\exp\left\{-4\int_0^s\|\bar X_r^{h,m}\|_\mV^2\dif r\right\}\|\bar X_s^{h,m-}\|^2_\mV\right)\\
&&+C_\t\sup\limits_{h\in{\bf D}^N}\sup\limits_{s\in[0,T]}\left(\exp\left\{-4\int_0^s\|\bar X_r^{h,m}\|_\mV^2\dif r\right\}|\bar X_s^{h,m-}|^2_\mH\right)\\
&\leq&C\t+C_\t\sup\limits_{h\in{\bf D}^N}\sup\limits_{s\in[0,T]}|\bar X_s^{h,m-}|^2_\mH,
\de
where we use (\ref{unifexpbarx}). As $m\rightarrow\infty$ first and then $\t\rightarrow 0$, (\ref{unifbarxhn4}) implies that
\ce
\lim\limits_{m\rightarrow\infty}\sup\limits_{h\in{\bf D}^N}\sup\limits_{s\in[0,T]}\left(\exp\left\{-2\int_0^s\|\bar X_r^{h,m}\|_\mV^2\dif r\right\}\|\bar X_s^{h,m-}\|_{L^\infty([0,1])}\right)=0.
\de
By the same deduction to that for the above limit, we obtain that 
\ce
\lim\limits_{n\rightarrow\infty}\sup\limits_{h\in{\bf D}^N}\sup\limits_{s\in[0,T]}\left(\exp\left\{-2\int_0^s\|\bar X_r^{h,n}\|_\mV^2\dif r\right\}\|\bar X_s^{h,n-}\|_{L^\infty([0,1])}\right)=0.
\de
The above deduction yields that
\be
\lim\limits_{n,m\rightarrow\infty}\sup\limits_{h\in{\bf D}^N}\left\{\sup\limits_{t\in[0,T]}\Sigma_t|\bar X_t^{h,n}-\bar{X}_t^{h,m}|_\mH^2+\int_0^T\Sigma_s\|\bar X_s^{h,n}-\bar{X}_s^{h,m}\|_\mV^2\dif s\right\}=0.
\label{sigmbarxhnbarxhm}
\ee
Besides, we notice that
\ce
\sup\limits_{h\in{\bf D}^N}\sup\limits_{t\in[0,T]}|\bar X_t^{h,n}-\bar{X}_t^{h,m}|_\mH^2&\leq&\(\sup\limits_{h\in{\bf D}^N}\sup\limits_{t\in[0,T]}\Sigma_t|\bar X_t^{h,n}-\bar{X}_t^{h,m}|_\mH^2\)\times\(\sup\limits_{h\in{\bf D}^N}\sup\limits_{t\in[0,T]}(\Sigma_t)^{-1}\)\\
&\leq&C\(\sup\limits_{h\in{\bf D}^N}\sup\limits_{t\in[0,T]}\Sigma_t|\bar X_t^{h,n}-\bar{X}_t^{h,m}|_\mH^2\),
\de
where (\ref{unifbarxhn}) is used. Similarly, it holds that
\ce
\sup\limits_{h\in{\bf D}^N}\int_0^T\|\bar X_s^{h,n}-\bar{X}_s^{h,m}\|_\mV^2\dif s\leq C\sup\limits_{h\in{\bf D}^N}\int_0^T\Sigma_s\|\bar X_s^{h,n}-\bar{X}_s^{h,m}\|_\mV^2\dif s,
\de
which together with (\ref{sigmbarxhnbarxhm}) implies (\ref{hnhmdiff}).
\end{proof}

For $\vartheta\in(0,1)$, let $W^{\vartheta, 2}([0, T], \mV^*)$ be the Sobolev space of all $v\in L^2([0, T], \mV^*)$ such that
$$
\int_0^T \int_0^T \frac{\|v(t)-v(s)\|_{\mV^*}^2}{|t-s|^{1+2\vartheta}}\dif t\dif s<\infty,
$$
endowed with the norm
$$
\|v\|_{W^{\vartheta, 2}([0, T], \mV^*)}^2:=\int_0^T\|v(t)\|_{\mV^*}^2\dif t+\int_0^T \int_0^T \frac{\|v(t)-v(s)\|_{\mV^*}^2}{|t-s|^{1+2\vartheta}}\dif t\dif s.
$$
The following lemma is from \cite[Lemma 4.3]{ntt}.

\bl\label{complemm}
Let $\Gamma$ be the space
\ce
\Gamma=L^\infty([0,T], \mH)\cap L^2([0,T], \mV)\cap W^{\vartheta, 2}([0, T], \mV^*)
\de
endowed with the natural norm. Then the embedding of $\Gamma$ in $L^2([0,T], \mH)$ is compact.
\el

Based on the above lemma, we obtain the following result.

\bl\label{barxhenbarxhn}
Suppose that $({\bf H}_{F_1, G_1})$, $({\bf H}^1_{F_2,G_2})$, $({\bf H}^2_{G_2})$ and $({\bf H}^3_{F_2, G_2})$ hold. Then for each $N<\infty$, any $\{h_\e\}\subset \mathbf{D}^N$ and $h\in \mathbf{D}^N$ satisfying $h_\e\to h$ as $\e\rightarrow0$,
\ce
\lim\limits_{\e\to0}\sup\limits_{t\in[0,T]}|\bar X_t^{h_\e,n}-\bar X_t^{h,n}|_\mH=0.
\de
\el
\begin{proof}
We divide the proof into two steps. In the first step, we prove that $\{\bar X^{h_\e,n}, \e>0\}$ is relatively compact in $L^2([0,T], \mH)$. Then the required limit is established in the second step.

{\bf Step 1.} We prove that $\{\bar X^{h_\e,n}, \e>0\}$ is relatively compact in $L^2([0,T], \mH)$.

In order to obtain that $\{\bar X^{h_\e,n}, \e>0\}$ is relatively compact in $L^2([0,T], \mH)$, by Lemma \ref{complemm} we only need to prove that $\{\bar X^{h_\e,n}, \e>0\}$ is bounded in $\Gamma$. By (\ref{unifbarxhn}), it holds that $\{\bar X^{h_\e,n}, \e>0\}$ is bounded in $L^\infty([0,T], \mH)\cap L^2([0,T], \mV)$. Then we show that $\{\bar X^{h_\e,n}, \e>0\}$ is bounded in $W^{\vartheta, 2}([0, T], \mV^*)$. Since $W^{1, 2}([0, T], \mV^*)\subset W^{\vartheta, 2}([0, T], \mV^*)$, where $W^{1, 2}([0, T], \mV^*)$ denotes the space of all $v\in L^2([0,T], \mV^*)$ such that $\frac{\dif v}{\dif t}\in L^2([0,T], \mV^*)$ with the norm $\|v\|_{W^{1, 2}([0, T], \mV^*)}:=\|v\|_{L^2([0,T], \mV^*)}+\|\frac{\dif v}{\dif t}\|_{L^2([0,T], \mV^*)}$, we prove that $\{\bar X^{h_\e,n}, \e>0\}$ is bounded in $W^{1, 2}([0, T], \mV^*)$.

Note that
\ce
\bar X_t^{h_\e,n}&=&x_0+\int_0^t A\bar X_s^{h_\e,n}\dif s+\int_0^t  B(\bar X_s^{h_\e,n},\bar X_s^{h_\e,n})\dif s+\int_0^t \bar F_1(\bar X_s^{h_\e,n})\dif s\\
&&+\int_0^t G_1(\bar X_s^{h_\e,n})\Pi_1h_\e(s)\dif s+\int_0^t n\bar X_s^{h_\e,n-}\dif s.
\de
Thus, it holds that
\ce
&&\|\bar X^{h_\e,n}\|^2_{W^{1, 2}([0, T], \mV^*)}\\
&\leq& \int_0^T\|\bar X_t^{h_\e,n}\|^2_{\mV^*}\dif t+5\int_0^T \|A\bar X_t^{h_\e,n}\|^2_{\mV^*}\dif t+5\int_0^T\|B(\bar X_t^{h_\e,n},\bar X_t^{h_\e,n})\|^2_{\mV^*}\dif t\\
&&+5\int_0^T\|\bar F_1(\bar X_t^{h_\e,n})\|^2_{\mV^*}\dif t+5\int_0^T\|G_1(\bar X_t^{h_\e,n})\Pi_1h_\e(t)\|^2_{\mV^*}\dif t+5\int_0^Tn^2\|\bar X_t^{h_\e,n-}\|^2_{\mV^*}\dif t\\
&=:&I_1+I_2+I_3+I_4+I_5+I_6.
\de
For $I_1$, by Lemma \ref{exisuniq} and (\ref{unifbarxhn}), we know that
\ce
I_1=\int_0^T|\bar X_t^{h_\e,n}|^2_{\mH}\dif t\leq \sup\limits_{t\in[0,T]}|\bar X_t^{h_\e,n}|^2_{\mH}T\leq C.
\de
For $I_2$ and $I_3$, the definition of $A$, Corollary \ref{Bbprop2} and (\ref{unifbarxhn}) imply that
\ce
&&I_2\leq C\int_0^T \|\bar X_t^{h_\e,n}\|^2_\mV\dif t\leq C,\\
&&I_3\leq 20\int_0^T |\bar X_t^{h_\e,n}|^2_{\mH}\|\bar X_t^{h_\e,n}\|^2_{\mV}\dif t\leq 20\sup\limits_{t\in[0,T]}|\bar X_t^{h_\e,n}|^2_{\mH}\int_0^T \|\bar X_t^{h_\e,n}\|^2_\mV\dif t\leq C.
\de
For $I_4$, by (\ref{barf1lip}) it holds that
\ce
I_4\leq C\int_0^T(1+|\bar X_t^{h_\e,n}|^2_{\mH})\dif t\leq C(1+\sup\limits_{t\in[0,T]}|\bar X_t^{h_\e,n}|^2_{\mH})T.
\de
For $I_5$, noticing that $\{h_\e\}\subset \mathbf{D}^N$, we infer that
\ce
I_5\leq 5\int_0^T\|G_1(\bar X_t^{h_\e,n})\|^2_{\cL_2(\mU_1,\mH)}|\Pi_1h_\e(t)|^2_{\mU_1}\dif t\leq C(1+\sup\limits_{t\in[0,T]}|\bar X_t^{h_\e,n}|^2_{\mH})\int_0^T|h_\e(t)|^2_{\sU}\dif t\leq CN.
\de
For $I_6$, it is easy to see that
\ce
I_6= 5\int_0^Tn^2|\bar X_t^{h_\e,n}|^2_{\mH}\dif t\leq 5 n^2\sup\limits_{t\in[0,T]}|\bar X_t^{h_\e,n}|^2_{\mH}T.
\de

Combining the above deduction, we conclude that $\{\bar X^{h_\e,n}, \e>0\}$ is bounded in $W^{1, 2}([0, T], \mV^*)$.

{\bf Step 2.} We prove that $\lim\limits_{\e\to0}\sup\limits_{t\in[0,T]}|\bar X_t^{h_\e,n}-\bar X_t^{h,n}|_\mH=0.$

First of all, define
\ce
\Sigma^n_t:=\exp\left\{-\b_2\int_0^t(\|\bar X_s^{h_\e,n}\|_\mV^2+\|\bar X_s^{h,n}\|_\mV^2)\dif s\right\}, \quad t\geq 0,
\de
where the constant $\b_2>0$ is determined later. By the chain rule we have that
\ce
\Sigma^n_t|\bar X_t^{h_\e,n}-\bar{X}_t^{h,n}|_\mH^2&=&-\b_2\int_0^t(\|\bar X_s^{h_\e,n}\|_\mV^2+\|\bar X_s^{h,n}\|_\mV^2)\Sigma^n_s|\bar X_s^{h_\e,n}-\bar{X}_s^{h,n}|_\mH^2\dif s\\
&&+2\int_0^t\Sigma^n_s{_{\mV^*}}\<A(\bar X_s^{h_\e,n}-\bar{X}_s^{h,n}),\bar X_s^{h_\e,n}-\bar{X}_s^{h,n}\>_\mV\dif s\\
&&+2\int_0^t\Sigma^n_s{_{\mV^*}}\<B(\bar X_s^{h_\e,n},\bar X_s^{h_\e,n})-B(\bar X_s^{h,n},\bar X_s^{h,n}),\bar X_s^{h_\e,n}-\bar{X}_s^{h,n}\>_\mV\dif s\\
&&+2\int_0^t\Sigma^n_s\<\bar F_1(\bar X_s^{h_\e,n})-\bar F_1(\bar X_s^{h,n}),\bar X_s^{h_\e,n}-\bar{X}_s^{h,n}\>_\mH\dif s\\
&&+2\int_0^t\Sigma^n_s\<\(G_1(\bar X_s^{h_\e,n})-G_1(\bar X_s^{h,n})\)\Pi_1h_\e(s),\bar X_s^{h_\e,n}-\bar{X}_s^{h,n}\>_\mH\dif s\\
&&+2\int_0^t\Sigma^n_s\<G_1(\bar X_s^{h,n})(\Pi_1h_\e(s)-\Pi_1h(s)),\bar X_s^{h_\e,n}-\bar{X}_s^{h,n}\>_\mH\dif s\\
&&+2\int_0^t\Sigma^n_s\<n\bar X_s^{h_\e,n-}-n\bar X_s^{h,n-},\bar X_s^{h_\e,n}-\bar{X}_s^{h,n}\>_\mH\dif s\\
&=:&J_1+J_2+J_3+J_4+J_5+J_6+J_7.
\de
By the similar deduction to that in Lemma \ref{unihdn}, it holds that
\ce
&&J_2=-2\int_0^t\Sigma^n_s|\bar X_s^{h_\e,n}-\bar{X}_s^{h,n}|^2_\mV\dif s,\\
&&J_3\leq \int_0^t\Sigma^n_s\|\bar X_s^{h_\e,n}-\bar{X}_s^{h,n}\|_\mV^2\dif s+8\int_0^t\Sigma^n_s|\bar X_s^{h_\e,n}-\bar{X}_s^{h,n}|_\mH^2(\|\bar X_s^{h_\e,n}\|^2_\mV+\|\bar X_s^{h,n}\|^2_\mV)\dif s,\\
&&J_4\leq 2C\int_0^t\Sigma^n_s|\bar X_s^{h_\e,n}-\bar{X}_s^{h,n}|^2_\mH\dif s,\\
&&J_5\leq \frac{1}{2}\sup\limits_{s\in[0,t]}\Sigma^n_s|\bar X_s^{h_\e,n}-\bar{X}_s^{h,n}|_\mH^2+2NL_1^2\int_0^t\Sigma^n_s|\bar X_s^{h_\e,n}-\bar{X}_s^{h,n}|_\mH^2\dif s.
\de
Besides, the fact that $|\bar X_t^{h_\e,n-}-\bar X_t^{h,n-}|_\mH\leq|\bar X_s^{h_\e,n}-\bar{X}_s^{h,n}|_\mH$ yields that
\ce
J_7\leq 2n\int_0^t\Sigma^n_s|\bar X_s^{h_\e,n}-\bar{X}_s^{h,n}|_\mH^2\dif s.
\de
Based on the above deduction, taking $\b_2\geq 8$, we obtain that
\ce
&&\sup\limits_{t\in[0,T]}\Sigma^n_t|\bar X_t^{h_\e,n}-\bar{X}_t^{h,n}|_\mH^2+\int_0^T\Sigma^n_s\|\bar X_s^{h_\e,n}-\bar{X}_s^{h,n}\|_\mV^2\dif s\\
&\leq&4(C+NL_1^2+n)\int_0^T\Sigma^n_s|\bar X_s^{h_\e,n}-\bar{X}_s^{h,n}|_\mH^2\dif s\\
&&+4\sup\limits_{t\in[0,T]}\left|\int_0^t\Sigma^n_s\<G_1(\bar X_s^{h,n})(\Pi_1h_\e(s)-\Pi_1h(s)),\bar X_s^{h_\e,n}-\bar{X}_s^{h,n}\>_\mH\dif s\right|,
\de
which together with the Gronwall inequality implies that
\ce
\sup\limits_{t\in[0,T]}\Sigma^n_t|\bar X_t^{h_\e,n}-\bar{X}_t^{h,n}|_\mH^2\leq C\sup\limits_{t\in[0,T]}\left|\int_0^t\Sigma^n_s\<G_1(\bar X_s^{h,n})(\Pi_1h_\e(s)-\Pi_1h(s)),\bar X_s^{h_\e,n}-\bar{X}_s^{h,n}\>_\mH\dif s\right|.
\de

Next, we prove that 
\ce
\lim\limits_{\e\to 0}\sup\limits_{t\in[0,T]}\left|\int_0^t\Sigma^n_s\<G_1(\bar X_s^{h,n})(\Pi_1h_\e(s)-\Pi_1h(s)),\bar X_s^{h_\e,n}-\bar{X}_s^{h,n}\>_\mH\dif s\right|=0.
\de
So, we only need to prove that for any sequence $\{\e_m, m\in\mN\}$ with $\lim\limits_{m\to\infty}\e_m=0$, 
\ce
\lim\limits_{m\to\infty}\sup\limits_{t\in[0,T]}\left|\int_0^t\Sigma^n_s\<G_1(\bar X_s^{h,n})(\Pi_1h_{\e_m}(s)-\Pi_1h(s)),\bar X_s^{h_{\e_m},n}-\bar{X}_s^{h,n}\>_\mH\dif s\right|=0.
\de

By {\bf Step 1}, we know that there exists a sequence $\{\bar X^{h_{\e_{m_i}},n}, i\in\mN\}$ and a $\check X^{h,n}\in L^2([0,T], \mH)$ such that $\bar X^{h_{\e_{m_i}},n}$ converges to $\check X^{h,n}$ in $L^2([0,T], \mH)$. So, it holds that
\ce
&&\sup\limits_{t\in[0,T]}\left|\int_0^t\Sigma^n_s\<G_1(\bar X_s^{h,n})(\Pi_1h_{\e_{m_i}}(s)-\Pi_1h(s)),\bar X_s^{h_{\e_{m_i}},n}-\bar{X}_s^{h,n}\>_\mH\dif s\right|\\
&\leq&\sup\limits_{t\in[0,T]}\left|\int_0^t\Sigma^n_s\<G_1(\bar X_s^{h,n})(\Pi_1h_{\e_{m_i}}(s)-\Pi_1h(s)),\bar X_s^{h_{\e_{m_i}},n}-\check{X}_s^{h,n}\>_\mH\dif s\right|\\
&&+\sup\limits_{t\in[0,T]}\left|\int_0^t\Sigma^n_s\<G_1(\bar X_s^{h,n})(\Pi_1h_{\e_{m_i}}(s)-\Pi_1h(s)),\check X_s^{h,n}-\bar{X}_s^{h,n}\>_\mH\dif s\right|\\
&=:&\cK_1+\cK_2.
\de
For $\cK_1$, the H\"older inequality and (\ref{f1g1linegrow}) yield that
\ce
\cK_1&\leq& \int_0^T|G_1(\bar X_s^{h,n})(\Pi_1h_{\e_{m_i}}(s)-\Pi_1h(s))|_\mH|\bar X_s^{h_{\e_{m_i}},n}-\check{X}_s^{h,n}|_\mH\dif s\\
&\leq&L_2^\prime(1+\sup\limits_{t\in[0,T]}|\bar X_s^{h,n}|_\mH)\(\int_0^T2(|h_{\e_{m_i}}(s)|^2_\sU+|h(s)|^2_\sU)\dif s\)^{1/2}\\
&&\times\( \int_0^T|\bar X_s^{h_{\e_{m_i}},n}-\check{X}_s^{h,n}|^2_\mH\dif s\)^{1/2}.
\de
From this, it follows that
\ce
\lim\limits_{i\to\infty}\cK_1=0.
\de
For $\cK_2$, noticing that $h_{\e_{m_i}}$ converges weakly to $h$ in $L^2([0,T],\sU)$, by the similar deduction to that for Lemma 5.3 in \cite{q}, we get that
\ce
\lim\limits_{i\to\infty}\cK_2=0.
\de
Finally, the above deduction implies the required limit. Thus, the proof is complete.
\end{proof}

Now, we define a mapping $\cG^0: C([0,T],\sU_0)\rightarrow C([0,T], \mH)$ by $\cG^0(\int_0^{\cdot}h(s)\dif s):=\bar{X}^{h}$ and establish the following result.

\bp\label{just1}
Suppose that $({\bf H}_{F_1, G_1})$, $({\bf H}^1_{F_2,G_2})$, $({\bf H}^2_{G_2})$ and $({\bf H}^3_{F_2, G_2})$ hold. Then for each $N<\infty$, any $\{h_\e\}\subset \mathbf{D}^N$ and $h\in \mathbf{D}^N$ satisfying $h_\e\to h$ as $\e\rightarrow0$, 
$$
\lim\limits_{\e\to0}\sup\limits_{t\in[0,T]}\left|\cG^0\left(\int_{0}^{\cdot}h_\e(s)\dif s\right)(t)-\cG^{0}\left(\int_{0}^{\cdot}h(s)\dif s\right)(t)\right|_\mH=0.
$$
\ep
\begin{proof}
Note that 
\ce
\cG^0\left(\int_{0}^{\cdot}h_\e(s)\dif s\right)=\bar X^{h_\e}, \quad \cG^{0}\left(\int_{0}^{\cdot}h(s)\dif s\right)=\bar X^{h}.
\de
Thus, it holds that
\ce
&&\sup\limits_{t\in[0,T]}\left|\cG^0\left(\int_{0}^{\cdot}h_\e(s)\dif s\right)(t)-\cG^{0}\left(\int_{0}^{\cdot}h(s)\dif s\right)(t)\right|_\mH=\sup\limits_{t\in[0,T]}\left|\bar X_t^{h_\e}-\bar X_t^{h}\right|_\mH\\
&\leq&\sup\limits_{t\in[0,T]}\left|\bar X_t^{h_\e}-\bar X_t^{h_\e,n}\right|_\mH+\sup\limits_{t\in[0,T]}\left|\bar X_t^{h_\e,n}-\bar X_t^{h,n}\right|_\mH+\sup\limits_{t\in[0,T]}\left|\bar X_t^{h,n}-\bar X_t^{h}\right|_\mH\\
&\leq&2\sup\limits_{h\in \mathbf{D}^N}\sup\limits_{t\in[0,T]}\left|\bar X_t^{h}-\bar X_t^{h,n}\right|_\mH+\sup\limits_{t\in[0,T]}\left|\bar X_t^{h_\e,n}-\bar X_t^{h,n}\right|_\mH.
\de
Letting $\e\to 0$ firstly and then $n\to\infty$, by Lemma \ref{barxhenbarxhn} and \ref{unihdn} we conclude that 
$$
\lim\limits_{\e\to0}\sup\limits_{t\in[0,T]}\left|\cG^0\left(\int_{0}^{\cdot}h_\e(s)\dif s\right)(t)-\cG^{0}\left(\int_{0}^{\cdot}h(s)\dif s\right)(t)\right|_\mH=0.
$$
The proof is complete.
\end{proof}

\subsection{Verification of Condition \ref{cond} $(ii)$}\label{vericond2}

In this subsection, we justify Condition \ref{cond} $(ii)$.

First of all, by the Yamada-Watanabe theorem there exists a measurable functional $\cG^\e: C([0,T], \sU_0)\to C([0,T], \mH)$ such that $X^{\e,\g}=\cG^\e(W)$. 

Consider the following controlled equation: for any $u_\e\in \cA_N$ 
\be\left\{\begin{array}{ll}
\dif X_t^{\e,\g,u_\e}=[AX_t^{\e,\g,u_\e}+B(X_t^{\e,\g,u_\e},X_t^{\e,\g,u_\e})+F_1(X_t^{\e,\g,u_\e},Y_t^{\e,\g,u_\e})]\dif t\\
\qquad\qquad\quad+\sqrt\e G_1(X_t^{\e,\g,u_\e})\dif W^1_t+G_1(X_t^{\e,\g,u_\e})\Pi_1 u_\e(t)\dif t+\dif K^{1,\e,\g,u_\e}_t,\\
\dif Y_t^{\e,\g,u_\e}=\frac{1}{\g}[AY_t^{\e,\g,u_\e}+F_2(X_t^{\e,\g,u_\e},Y_t^{\e,\g,u_\e})]\dif t+\frac{\sqrt \e}{\sqrt{\g\e}} G_2(X_t^{\e,\g,u_\e},Y_t^{\e,\g,u_\e})\dif W^2_t\\
\qquad\qquad\quad+\frac{1}{\sqrt{\g\e}}G_2(X_t^{\e,\g,u_\e},Y_t^{\e,\g,u_\e})\Pi_2 u_\e(t)\dif t+\dif K^{2,\e,\g,u_\e}_t,\\
X_t^{\e,\g,u_\e}(\xi)\geq 0, \quad Y_t^{\e,\g,u_\e}(\xi)\geq 0,\quad \xi\in[0,1],\\
X_t^{\e,\g,u_\e}(0)=X_t^{\e,\g,u_\e}(1)=Y_t^{\e,\g,u_\e}(0)=Y_t^{\e,\g,u_\e}(1)=0,\\
X_0^{\e,\g,u_\e}=x_0, \quad Y_0^{\e,\g,u_\e}=y_0.
\end{array}
\right.
\label{contequa}
\ee
By the Girsanov theorem and Theorem \ref{wellpose}, we know that the system (\ref{contequa}) has a unique solution $(X^{\e,\g,u_\e}, K^{1,\e,\g,u_\e}, Y^{\e,\g,u_\e}, K^{2,\e,\g,u_\e})$. Moreover, 
$$
X^{\e,\g,u_\e}=\cG^\e\left(W(\cdot)+\frac{1}{\sqrt{\varepsilon}} \int_0^{\cdot} u_\e(s)\dif s\right).
$$
In order to justify Condition \ref{cond} $(ii)$, we prepare some estimates.

\bl
Under $({\bf H}_{F_1, G_1})$, $({\bf H}^1_{F_2,G_2})$, $({\bf H}^3_{F_2,G_2})$ and $({\bf H}^4_{G_2})$ it holds that
\be
&&\mE\left(\sup _{t \in[0, T]}\left|X_t^{\e,\g,u_\e}\right|_\mH^2\right)+\mE \int_0^T\left\|X_t^{\e,\g,u_\e}\right\|_\mV^2\dif t \leq C\left(1+|x_0|_\mH^2+|y_0|_\mH^2\right),\label{xeguees}\\
&&\mE\int_0^T\left|Y_t^{\e,\g,u_\e}\right|_\mH^2\dif t \leq C\left(1+|x_0|_\mH^2+|y_0|_\mH^2\right),\label{yeguees}
\ee
where the constant $C>0$ is independent of $\e,\g$.
\el
\begin{proof}
By the It\^o formula, it holds that for any $t\in[0,T]$
\ce
\left|X_t^{\e,\g,u_\e}\right|_\mH^2&=&|x_0|_\mH^2+2\int_0^t {_{\mV^*}}\<AX_s^{\e,\g,u_\e},X_s^{\e,\g,u_\e}\>_\mV \dif s+2 \int_0^t{_{\mV^*}}\<B(X_s^{\e,\g,u_\e},X_s^{\e,\g,u_\e}),X_s^{\e,\g,u_\e}\>_\mV\dif s\\
&& +2 \int_0^t\<F_1(X_s^{\e,\g,u_\e},Y_s^{\e,\g,u_\e}), X_s^{\e,\g,u_\e}\>_\mH\dif s+2 \sqrt{\e} \int_0^t\<X_s^{\e,\g,u_\e}, G_1(X_s^{\e,\g,u_\e})\dif W^1_s\>_\mH\\
&& +2 \int_0^t\<X_s^{\e,\g,u_\e}, G_1(X_s^{\e,\g,u_\e})\Pi_1 u_\e(s)\>_\mH\dif s+2 \int_0^t\int_0^1X_s^{\e,\g,u_\e}(\xi)K^{1,\e,\g,u_\e}(\dif s,\dif \xi)\\
&&+\e \int_0^t\|G_1(X_s^{\e,\g,u_\e})\|_{\cL_2(\mU_1,\mH)}^2\dif s. 
\de
Note that 
\ce
&&2{_{\mV^*}}\<AX_s^{\e,\g,u_\e},X_s^{\e,\g,u_\e}\>_\mV=-2\|X_s^{\e,\g,u_\e}\|^2_\mV,\\
&&2{_{\mV^*}}\<B(X_s^{\e,\g,u_\e},X_s^{\e,\g,u_\e}),X_s^{\e,\g,u_\e}\>_\mV=0,\\
&&2 \int_0^t\int_0^1X_s^{\e,\g,u_\e}(\xi)K^{1,\e,\g,u_\e}(\dif s,\dif \xi)=0,
\de
where we use the definition of $A$, Corollary \ref{Bbprop2} and Definition \ref{soludefi} $(iv)$. Thus, by (\ref{f1g1linegrow}), the H\"older inequality and $u_\e\in \cA_N$ we infer that
\ce
\left|X_t^{\e,\g,u_\e}\right|_\mH^2+2\int_0^t\|X_s^{\e,\g,u_\e}\|_\mV^2\dif s&\leq& |x_0|_\mH^2+C+C\int_0^t|X_s^{\e,\g,u_\e}|_\mH^2\dif s+C\int_0^t|Y_s^{\e,\g,u_\e}|_\mH^2\dif s\\
 &&+2 \sqrt{\e} \left|\int_0^t\<X_s^{\e,\g,u_\e}, G_1(X_s^{\e,\g,u_\e})\dif W^1_s\>_\mH\right|+\frac{1}{4}\sup\limits_{s\in[0,t]}|X_s^{\e,\g,u_\e}|_\mH^2\\
 &&+(4N+1)\int_0^t\|G_1(X_s^{\e,\g,u_\e})\|_{\cL_2(\mU_1,\mH)}^2\dif s. 
\de
The Burkholder-Davis-Gundy inequality implies that
\ce
&&\mE\sup\limits_{s\in[0,t]}|X_s^{\e,\g,u_\e}|_\mH^2+2\mE\int_0^t\|X_s^{\e,\g,u_\e}\|_\mV^2\dif s\no\\
&\leq& |x_0|_\mH^2+C+C\mE\int_0^t|X_s^{\e,\g,u_\e}|_\mH^2\dif s+C\mE\int_0^t|Y_s^{\e,\g,u_\e}|_\mH^2\dif s\no\\
&&+2 \sqrt{\e}C \mE\left(\int_0^t|X_s^{\e,\g,u_\e}|_\mH^2\|G_1(X_s^{\e,\g,u_\e})\|_{\cL_2(\mU_1,\mH)}^2\dif s\right)^{1/2}+\frac{1}{4}\mE\sup\limits_{s\in[0,t]}|X_s^{\e,\g,u_\e}|_\mH^2\no\\
&&+(4N+1)\mE\int_0^t\|G_1(X_s^{\e,\g,u_\e})\|_{\cL_2(\mU_1,\mH)}^2\dif s\no\\
&\leq&|x_0|_\mH^2+C+C\mE\int_0^t|X_s^{\e,\g,u_\e}|_\mH^2\dif s+C\mE\int_0^t|Y_s^{\e,\g,u_\e}|_\mH^2\dif s\no\\
&&+\frac{1}{2}\mE\sup\limits_{s\in[0,t]}|X_s^{\e,\g,u_\e}|_\mH^2+C\mE\int_0^t\|G_1(X_s^{\e,\g,u_\e})\|_{\cL_2(\mU_1,\mH)}^2\dif s.
\de
and furthermore by (\ref{f1g1linegrow})
\be
&&\mE\sup\limits_{s\in[0,t]}|X_s^{\e,\g,u_\e}|_\mH^2+\mE\int_0^t\|X_s^{\e,\g,u_\e}\|_\mV^2\dif s\no\\
&\leq& C(|x_0|_\mH^2+1)+C\mE\int_0^t|X_s^{\e,\g,u_\e}|_\mH^2\dif s+C\mE\int_0^t|Y_s^{\e,\g,u_\e}|_\mH^2\dif s.
\label{xegueyegue1}
\ee

Next, we compute $\mE\int_0^t|Y_s^{\e,\g,u_\e}|_\mH^2\dif s$. From the It\^o formula, it follows that for any $l_3>0$
\ce
e^{l_3 t}|Y_t^{\e,\g,u_\e}|_\mH^2&=&|y_0|_\mH^2+l_3\int_0^te^{l_3 s}|Y_s^{\e,\g,u_\e}|_\mH^2\dif s+\frac{2}{\g}\int_0^te^{l_3 s}{_{\mV^*}}\<AY_s^{\e,\g,u_\e},Y_s^{\e,\g,u_\e}\>_\mV \dif s\\
&&+\frac{2}{\g}\int_0^te^{l_3 s}\<F_2(X_s^{\e,\g,u_\e},Y_s^{\e,\g,u_\e}), Y_s^{\e,\g,u_\e}\>_\mH\dif s\\
&&+\frac{2\sqrt \e}{\sqrt{\g\e}}\int_0^te^{l_3 s}\<Y_s^{\e,\g,u_\e},G_2(X_s^{\e,\g,u_\e},Y_s^{\e,\g,u_\e})\dif W^2_s\>_\mH\\
&&+\frac{2}{\sqrt{\g\e}}\int_0^te^{l_3 s}\<Y_s^{\e,\g,u_\e},G_2(X_s^{\e,\g,u_\e},Y_s^{\e,\g,u_\e})\Pi_2 u_\e(s)\>_\mH\dif s\\
&&+2\int_0^t\int_0^1e^{l_3 s}Y_s^{\e,\g,u_\e}(\xi)K^{2,\e,\g,u_\e}(\dif s,\dif\xi)\\
&&+\frac{1}{\g}\int_0^te^{l_3 s}\|G_2(X_s^{\e,\g,u_\e},Y_s^{\e,\g,u_\e})\|_{\cL_2(\mU_2,\mH)}^2\dif s.
\de
Besides, (\ref{invameasex}) implies that
\ce
&&\frac{2}{\g}{_{\mV^*}}\<AY_s^{\e,\g,u_\e},Y_s^{\e,\g,u_\e}\>_\mV+\frac{2}{\g}\<F_2(X_s^{\e,\g,u_\e},Y_s^{\e,\g,u_\e}), Y_s^{\e,\g,u_\e}\>_\mH+\frac{1}{\g}\|G_2(X_s^{\e,\g,u_\e},Y_s^{\e,\g,u_\e})\|_{\cL_2(\mU_2,\mH)}^2\\
&\leq&-\frac{2\a}{\g}|Y_s^{\e,\g,u_\e}|_\mH^2+\frac{C}{\g}(1+|X_s^{\e,\g,u_\e}|_\mH^2).
\de
And the Young inequality yields that
\ce
&&\frac{2}{\sqrt{\g\e}}\<Y_s^{\e,\g,u_\e},G_2(X_s^{\e,\g,u_\e},Y_s^{\e,\g,u_\e})\Pi_2 u_\e(s)\>_\mH\\
&\leq& \frac{\a}{\g}|Y_s^{\e,\g,u_\e}|_\mH^2+\frac{C}{\e}\|G_2(X_s^{\e,\g,u_\e},Y_s^{\e,\g,u_\e})\|_{\cL_2(\mU_2,\mH)}^2|\Pi_2 u_\e(s)|_{\mU_2}^2.
\de
By Definition \ref{soludefi} $(iv)$, it holds that
\ce
2\int_0^t\int_0^1e^{l_3 s}Y_s^{\e,\g,u_\e}(\xi)K^{2,\e,\g,u_\e}(\dif s,\dif\xi)=0.
\de
So, we infer that
\ce
\mE e^{l_3 t}|Y_t^{\e,\g,u_\e}|_\mH^2&\leq& |y_0|_\mH^2+(l_3-\frac{2\a}{\g}+ \frac{\a}{\g})\mE\int_0^te^{l_3 s}|Y_s^{\e,\g,u_\e}|_\mH^2\dif s\\
&&+\frac{C}{\g}\mE\int_0^te^{l_3 s}(1+|X_s^{\e,\g,u_\e}|_\mH^2)\dif s\\
&&+\frac{C}{\e}\mE\int_0^te^{l_3 s}\|G_2(X_s^{\e,\g,u_\e},Y_s^{\e,\g,u_\e})\|_{\cL_2(\mU_2,\mH)}^2|\Pi_2 u_\e(s)|_{\mU_2}^2\dif s\\
&\leq& |y_0|_\mH^2+(l_3-\frac{2\a}{\g}+ \frac{\a}{\g})\mE\int_0^te^{l_3 s}|Y_s^{\e,\g,u_\e}|_\mH^2\dif s\\
&&+\frac{C}{\g}(1+\mE\sup\limits_{s\in[0,t]}|X_s^{\e,\g,u_\e}|_\mH^2)\int_0^te^{l_3 s}\dif s\\
&&+\frac{C}{\e}\mE\int_0^te^{l_3 s}L^2_8(1+|X_s^{\e,\g,u_\e}|_\mH)^2|\Pi_2 u_\e(s)|_{\mU_2}^2\dif s.
\de
By taking $l_3= \frac{\a}{\g}$, it holds that
\ce
\mE|Y_t^{\e,\g,u_\e}|_\mH^2&\leq& |y_0|_\mH^2+\frac{C}{\a}(1+\mE\sup\limits_{s\in[0,t]}|X_s^{\e,\g,u_\e}|_\mH^2)\\
&&+\frac{C}{\e}L^2_8\mE(1+\sup\limits_{s\in[0,T]}|X_s^{\e,\g,u_\e}|_\mH)^2\int_0^te^{-\frac{\a}{\g}(t-s)}|\Pi_2 u_\e(s)|_{\mU_2}^2\dif s.
\de
Moreover, integrating on both sides, we conclude that
\be
\int_0^t\mE|Y_r^{\e,\g,u_\e}|_\mH^2\dif r&\leq& C(|y_0|_\mH^2+1)+\frac{C}{\a}\int_0^t\mE\sup\limits_{s\in[0,r]}|X_s^{\e,\g,u_\e}|_\mH^2\dif r\no\\
&&+\frac{C}{\e}L^2_8\mE(1+\sup\limits_{s\in[0,T]}|X_s^{\e,\g,u_\e}|_\mH)^2\int_0^t\int_0^re^{-\frac{\a}{\g}(r-s)}|\Pi_2 u_\e(s)|_{\mU_2}^2\dif s\dif r\no\\
&\leq&C(|y_0|_\mH^2+1)+\frac{C}{\a}\int_0^t\mE\sup\limits_{s\in[0,r]}|X_s^{\e,\g,u_\e}|_\mH^2\dif r\no\\
&&+C(\frac{\g}{\e})L^2_8\mE(1+\sup\limits_{s\in[0,T]}|X_s^{\e,\g,u_\e}|_\mH)^2\int_0^t|u_\e(s)|_\sU^2\dif s\no\\
&\leq&C(|y_0|_\mH^2+1)+\frac{C}{\a}\int_0^t\mE\sup\limits_{s\in[0,r]}|X_s^{\e,\g,u_\e}|_\mH^2\dif r\no\\
&&+CN(\frac{\g}{\e})L^2_8\mE(1+\sup\limits_{s\in[0,T]}|X_s^{\e,\g,u_\e}|_\mH)^2,
\label{yegue}
\ee
where we make use of $u_\e\in \cA_N$.

Finally, by inserting (\ref{yegue}) into (\ref{xegueyegue1})  and noticing $\lim\limits_{\e\to 0}\frac{\g}{\e}=0$, the Gronwall inequality implies (\ref{xeguees}). Then by (\ref{xeguees}) and (\ref{yegue}), we obtain (\ref{yeguees}). The proof is complete.
\end{proof}

Next, in order to give the estimate about the path of $X^{\e,\g,u_\e}$, for any large $R>0$ we define the stopping time 
\ce
\tau_1:=\inf\{t>0, |X_t^{\e,\g,u_\e}|_\mH>R\},
\de
and have the following result.

\bl
Under $({\bf H}_{F_1, G_1})$, $({\bf H}^1_{F_2,G_2})$, $({\bf H}^3_{F_2,G_2})$ and $({\bf H}^4_{G_2})$ it holds that for any $\d>0$ small enough,
\be
\mE\left[\int_0^{T\wedge\tau_1}|X_t^{\e,\g,u_\e}-X_{t(\d)}^{\e,\g,u_\e}|_\mH^2\dif t\right]\leq C_{R}\d^{1/2}\left(1+|x_0|_\mH^2+|y_0|_\mH^2\right),
\label{xeguettd}
\ee
where $t(\d):=[\frac{t}{\d}]\d$, $[\frac{t}{\d}]$ denotes the largest integer which is less than $\frac{t}{\d}$ and the constant $C_{R}>0$ depends on $R$.
\el
\begin{proof} 
First of all, we notice that
\ce
&&\mE\left[\int_0^{T\wedge\tau_1}|X_t^{\e,\g,u_\e}-X_{t(\d)}^{\e,\g,u_\e}|_\mH^2\dif t\right]=\mE\left[\int_0^{T}|X_t^{\e,\g,u_\e}-X_{t(\d)}^{\e,\g,u_\e}|_\mH^2 I_{t\leq\tau_1}\dif t\right]\\
&=& \mE\left[\int_0^\d|X_t^{\e,\g,u_\e}-X_{t(\d)}^{\e,\g,u_\e}|_\mH^2 I_{t\leq\tau_1}\dif t\right]+\mE\left[\int_\d^{T}|X_t^{\e,\g,u_\e}-X_{t(\d)}^{\e,\g,u_\e}|_\mH^2 I_{t\leq\tau_1}\dif t\right]\\
&\leq&\mE\left[\int_0^\d|X_t^{\e,\g,u_\e}-x_0|_\mH^2\dif t\right]+2\mE\left[\int_\d^{T}|X_t^{\e,\g,u_\e}-X_{t-\d}^{\e,\g,u_\e}|_\mH^2 I_{t\leq\tau_1}\dif t\right]\\
&&+2\mE\left[\int_\d^{T}|X_{t-\d}^{\e,\g,u_\e}-X_{t(\d)}^{\e,\g,u_\e}|_\mH^2 I_{t\leq\tau_1}\dif t\right]\\
&\leq&C\left(1+|x_0|_\mH^2+|y_0|_\mH^2\right)\d+2\mE\left[\int_\d^{T}|X_t^{\e,\g,u_\e}-X_{t-\d}^{\e,\g,u_\e}|_\mH^2 I_{t\leq\tau_1}\dif t\right]\\
&&+2\mE\left[\int_\d^{T}|X_{t-\d}^{\e,\g,u_\e}-X_{t(\d)}^{\e,\g,u_\e}|_\mH^2 I_{t\leq\tau_1}\dif t\right],
\de
where we use (\ref{xeguees}). Since the estimates for the second term and the third term in the right hand of the above inequality are similar, we only deal with the second term.

Next, we compute $\mE\left[\int_\d^{T}|X_t^{\e,\g,u_\e}-X_{t-\d}^{\e,\g,u_\e}|_\mH^2 I_{t\leq\tau_1}\dif t\right]$. Applying the It\^o formula to $|X_s^{\e,\g,u_\e}-X_{t-\d}^{\e,\g,u_\e}|_\mH^2$ for $s\in[t-\d,t]$, we have that
\ce
&&|X_t^{\e,\g,u_\e}-X_{t-\d}^{\e,\g,u_\e}|_\mH^2\\
&=&2\int_{t-\d}^t{_{\mV^*}}\<A X_s^{\e,\g,u_\e}, X_s^{\e,\g,u_\e}-X_{t-\d}^{\e,\g,u_\e}\>_\mV\dif s\\
&&+2\int_{t-\d}^t{_{\mV^*}}\<B(X_s^{\e,\g,u_\e},X_s^{\e,\g,u_\e}), X_s^{\e,\g,u_\e}-X_{t-\d}^{\e,\g,u_\e}\>_\mV\dif s\\
&& +2 \int_{t-\d}^t\<F_1(X_s^{\e,\g,u_\e},Y_s^{\e,\g,u_\e}), X_s^{\e,\g,u_\e}-X_{t-\d}^{\e,\g,u_\e}\>_\mH\dif s \\
&& +2 \sqrt{\e} \int_{t-\d}^t\<X_s^{\e,\g,u_\e}-X_{t-\d}^{\e,\g,u_\e}, G_1(X_s^{\e,\g,u_\e})\dif W^1_s\>_\mH \\
&& +2 \int_{t-\d}^t\<X_s^{\e,\g,u_\e}-X_{t-\d}^{\e,\g,u_\e}, G_1(X_s^{\e,\g,u_\e})\Pi_1 u_\e(s)\>_\mH\dif s\\
&&+2 \int_{t-\d}^t\int_0^1(X_s^{\e,\g,u_\e}(\xi)-X_{t-\d}^{\e,\g,u_\e}(\xi))K^{1,\e,\g,u_\e}(\dif s,\dif \xi)\\
&&+\e \int_{t-\d}^t\|G_1(X_s^{\e,\g,u_\e})\|_{\cL_2(\mU_1,\mH)}^2\dif s.
\de
By the definition of $A$, it holds that
\ce
&&2\int_{t-\d}^t{_{\mV^*}}\<A X_s^{\e,\g,u_\e}, X_s^{\e,\g,u_\e}-X_{t-\d}^{\e,\g,u_\e}\>_\mV\dif s\\
&\leq&-2\int_{t-\d}^t\|X_s^{\e,\g,u_\e}\|_\mV^2\dif s+2\int_{t-\d}^t\|X_s^{\e,\g,u_\e}\|_\mV\|X_{t-\d}^{\e,\g,u_\e}\|_\mV\dif s\\
&\leq&2\int_{t-\d}^t\|X_s^{\e,\g,u_\e}\|_\mV\cdot\|X_{t-\d}^{\e,\g,u_\e}\|_\mV\dif s.
\de
And Corollary \ref{Bbprop2} implies that
\ce
&&2\int_{t-\d}^t{_{\mV^*}}\<B(X_s^{\e,\g,u_\e},X_s^{\e,\g,u_\e}), X_s^{\e,\g,u_\e}-X_{t-\d}^{\e,\g,u_\e}\>_\mV\dif s\\
&=&-2\int_{t-\d}^t{_{\mV^*}}\<B(X_s^{\e,\g,u_\e},X_s^{\e,\g,u_\e}), X_{t-\d}^{\e,\g,u_\e}\>_\mV\dif s\\
&\leq&4\int_{t-\d}^t\|X_s^{\e,\g,u_\e}\|_\mV\cdot|X_s^{\e,\g,u_\e}|_\mH\cdot\|X_{t-\d}^{\e,\g,u_\e}\|_\mV\dif s.
\de
From (\ref{f1g1linegrow}), it follows that
\ce
&&2 \int_{t-\d}^t\<F_1(X_s^{\e,\g,u_\e},Y_s^{\e,\g,u_\e}), X_s^{\e,\g,u_\e}-X_{t-\d}^{\e,\g,u_\e}\>_\mH\dif s \\
&\leq&2L'_1\int_{t-\d}^t(1+|X_s^{\e,\g,u_\e}|_\mH+|Y_s^{\e,\g,u_\e}|_\mH)(|X_s^{\e,\g,u_\e}|_\mH+|X_{t-\d}^{\e,\g,u_\e}|_\mH)\dif s\\
&\leq&3L'_1\d+7L'_1\d\sup\limits_{s\in[0,T]}|X_s^{\e,\g,u_\e}|_\mH^2+3L'_1\int_{t-\d}^t|Y_s^{\e,\g,u_\e}|_\mH^2\dif s,
\de
and
\ce
\e \int_{t-\d}^t\|G_1(X_s^{\e,\g,u_\e})\|_{\cL_2(\mU_1,\mH)}^2\dif s&\leq& 2\e L^{\prime 2}_2\int_{t-\d}^t(1+|X_s^{\e,\g,u_\e}|_\mH^2)\dif s\\
&\leq& 2\e L^{\prime 2}_2\d+2\e L^{\prime 2}_2\d\sup\limits_{s\in[0,T]}|X_s^{\e,\g,u_\e}|_\mH^2.
\de
Definition \ref{soludefi} yields that
\ce
2 \int_{t-\d}^t\int_0^1(X_s^{\e,\g,u_\e}(\xi)-X_{t-\d}^{\e,\g,u_\e}(\xi))K^{1,\e,\g,u_\e}(\dif s,\dif \xi)\leq 0.
\de

So, we infer that
\ce
|X_t^{\e,\g,u_\e}-X_{t-\d}^{\e,\g,u_\e}|_\mH^2&\leq& 2\int_{t-\d}^t\|X_s^{\e,\g,u_\e}\|_\mV(1+2|X_s^{\e,\g,u_\e}|_\mH)\|X_{t-\d}^{\e,\g,u_\e}\|_\mV\dif s\\
&&+\left(C\d+C\d\sup\limits_{s\in[0,T]}|X_s^{\e,\g,u_\e}|_\mH^2+3L'_1\int_{t-\d}^t|Y_s^{\e,\g,u_\e}|_\mH^2\dif s\right)\\
&& +2 \sqrt{\e} \left|\int_{t-\d}^t\<X_s^{\e,\g,u_\e}-X_{t-\d}^{\e,\g,u_\e}, G_1(X_s^{\e,\g,u_\e})\dif W^1_s\>_\mH\right| \\
&& +2 \int_{t-\d}^t|X_s^{\e,\g,u_\e}-X_{t-\d}^{\e,\g,u_\e}|_\mH\cdot|G_1(X_s^{\e,\g,u_\e})\Pi_1 u_\e(s)|_\mH\dif s\\
&=:& I_1(t)+I_2(t)+I_3(t)+I_4(t).
\de

For $I_1(t)$, by the H\"older inequality, it holds that
\ce
\mE\left[\int_\d^{T}I_1(t)I_{t\leq\tau_1}\dif t\right]&\leq&2(1+2R)\left[\mE\int_\d^{T}\int_{t-\d}^t\|X_s^{\e,\g,u_\e}\|_\mV^2\dif s\dif t\right]^{1/2}\\
&&\qquad\qquad \times\left[\mE\int_\d^{T}\int_{t-\d}^t\|X_{t-\d}^{\e,\g,u_\e}\|_\mV^2\dif s\dif t\right]^{1/2}\\
&\leq&2(1+2R)\d \mE\int_0^{T}\|X_s^{\e,\g,u_\e}\|_\mV^2\dif s\\
&\leq&2(1+2R)\d C\left(1+|x_0|_\mH^2+|y_0|_\mH^2\right),
\de
where we use the Fubini theorem and (\ref{xeguees}) in the second and third inequalities, respectively.

For $I_2(t)$, by the Fubini theorem, (\ref{xeguees}) and (\ref{yeguees}), we have that
\ce
\mE\left[\int_\d^{T}I_2(t)I_{t\leq\tau_1}\dif t\right]&\leq&CT\d+CT\d\mE\sup\limits_{s\in[0,T]}|X_s^{\e,\g,u_\e}|_\mH^2+3L'_1\mE\int_\d^{T}\int_{t-\d}^t|Y_s^{\e,\g,u_\e}|_\mH^2\dif s\dif t\\
&\leq&CT\d+CT\d\mE\sup\limits_{s\in[0,T]}|X_s^{\e,\g,u_\e}|_\mH^2+3L'_1\d\mE\int_0^{T}|Y_s^{\e,\g,u_\e}|_\mH^2\dif s\\
&\leq&C\d\left(1+|x_0|_\mH^2+|y_0|_\mH^2\right).
\de

For $I_3(t)$, by the Burkholder-Davis-Gundy inequality, the H\"older inequality and (\ref{f1g1linegrow}), it holds that
\ce
&&\mE\left[\int_\d^{T}I_3(t)I_{t\leq\tau_1}\dif t\right]\\
&\leq&2 \sqrt{\e}C\int_\d^{T}\mE\left(\int_{t-\d}^t|X_s^{\e,\g,u_\e}-X_{t-\d}^{\e,\g,u_\e}|_\mH^2\|G_1(X_s^{\e,\g,u_\e})\|_{\cL_2(\mU_1,\mH)}^2I_{t\leq\tau_1}\dif s\right)^{1/2}\dif t\\
&\leq&2 \sqrt{\e}CT^{1/2}\left[\mE\int_\d^{T}\int_{t-\d}^t(|X_s^{\e,\g,u_\e}|_\mH^2+|X_{t-\d}^{\e,\g,u_\e}|_\mH^2)(1+|X_s^{\e,\g,u_\e}|_\mH^2)I_{t\leq\tau_1}\dif s\dif t\right]^{1/2}\\
&\leq&2 \sqrt{\e}CT^{1/2}(1+R^2)^{1/2}\left[\mE\int_\d^{T}\int_{t-\d}^t(|X_s^{\e,\g,u_\e}|_\mH^2+|X_{t-\d}^{\e,\g,u_\e}|_\mH^2)\dif s\dif t\right]^{1/2}\\
&\leq&2 \sqrt{\e}CT(1+R^2)^{1/2}\d^{1/2}\left[\mE\sup\limits_{s\in[0,T]}|X_s^{\e,\g,u_\e}|_\mH^2\right]^{1/2}\\
&\leq&2 \sqrt{\e}CT(1+R^2)^{1/2}\d^{1/2}\left(1+|x_0|_\mH^2+|y_0|_\mH^2\right)^{1/2}.
\de

For $I_4(t)$, (\ref{f1g1linegrow}) implies that
\ce
&&\mE\left[\int_\d^{T}I_4(t)I_{t\leq\tau_1}\dif t\right]\\
&\leq&2L'_2\mE\int_{\d}^T \int_{t-\d}^t(|X_s^{\e,\g,u_\e}|_\mH+|X_{t-\d}^{\e,\g,u_\e}|_\mH)(1+|X_s^{\e,\g,u_\e}|_\mH)|\Pi_1 u_\e(s)|_{\mU_1}\dif s I_{t \leq \tau_1}\dif t\\
&\leq&2L'_2(1+R)\left(\mE\int_{\d}^T \int_{t-\d}^t(|X_s^{\e,\g,u_\e}|_\mH+|X_{t-\d}^{\e,\g,u_\e}|_\mH)^2\dif s\dif t\right)^{1/2}\left(\mE\int_{\d}^T \int_{t-\d}^t|\Pi_1 u_\e(s)|_{\mU_1}^2\dif s\dif t\right)^{1/2}\\
&\leq&4L'_2(1+R)T^{1/2}\d\left[\mE\sup\limits_{s\in[0,T]}|X_s^{\e,\g,u_\e}|_\mH^2\right]^{1/2}\left(\mE\int_0^T|u_\e(s)|_{\sU}^2\dif s\right)^{1/2}\\
&\leq&4L'_2(1+R)T^{1/2}\d CN^{1/2}\left(1+|x_0|_\mH^2+|y_0|_\mH^2\right)^{1/2}.
\de

Finally, collecting the above deduction and using $z^{1/2}\leq 1+z$ for $z\geq 0$, we obtain that
\ce
2\mE\left[\int_\d^{T}|X_{t-\d}^{\e,\g,u_\e}-X_{t(\d)}^{\e,\g,u_\e}|_\mH^2 I_{t\leq\tau_1}\dif t\right]\leq C_{R}\d^{1/2}\left(1+|x_0|_\mH^2+|y_0|_\mH^2\right),
\de
which completes the proof.
\end{proof}

\bl
Under $({\bf H}_{F_1, G_1})$, $({\bf H}^1_{F_2,G_2})$, $({\bf H}^2_{G_2})$ and $({\bf H}^3_{F_2, G_2})$ it holds that for any $u\in\cA_N$ and $\d>0$ small enough,
\be
&&\left(\sup _{t \in[0, T]}\left|\bar X_t^{u}\right|_\mH^2\right)+\int_0^T\left\|\bar X_t^{u}\right\|_\mV^2\dif t\leq C\left(1+|x_0|_\mH^2\right), a.s.,\label{barxue}\\
&&\int_0^{T}|\bar X_t^{u}-\bar X_{t(\d)}^{u}|_\mH^2\dif t\leq C\d\left(1+|x_0|_\mH^2\right), a.s.,\label{barxueincr}
\ee
where $\bar X^{u}$ is the solution of Eq.(\ref{e0conteq}) with replacing $h$ by $u$ and $C>0$ is independent of $\omega$.
\el

Since the proofs of (\ref{barxue}) and (\ref{barxueincr}) are similar to that for (\ref{xeguees}) and (\ref{xeguettd}), respectively, we omit them.

Finally we introduce the following auxiliary SPDE with reflection: 
\be\left\{\begin{array}{l}
\dif\hat{Y}_{t}^{\e,\g,u_\e}=\frac{1}{\g}[A\hat{Y}_{t}^{\e,\g,u_\e}+F_{2}(X_{t(\d)}^{\e,\g,u_\e},\hat{Y}_{t}^{\e,\g,u_\e})]\dif t+\frac{1}{\sqrt{\g}}G_{2}(X_{t(\d)}^{\e,\g,u_\e},\hat{Y}_{t}^{\e,\g,u_\e})\dif W^2_{t}\\
\qquad\qquad\qquad+\dif \hat{K}_{t}^{2,\e,\g,u_\e}, \\
\hat{Y}_{t}^{\e,\g,u_\e}(\xi)\geq 0,\quad \xi\in[0,1],\\
\hat{Y}_{t}^{\e,\g,u_\e}(0)=\hat{Y}_{t}^{\e,\g,u_\e}(1)=0,\\
\hat{Y}_{0}^{\e,\g,u_\e}=y_0.
\end{array}
\right.
\label{hatzu}
\ee
Under $({\bf H}^1_{F_2,G_2})$, the above SPDE with reflection has a unique solution $(\hat{Y}_{t}^{\e,\g,u_\e}, \hat{K}_{t}^{\e,\g,u_\e})$ (\cite{xz}). 

\bl
Under $({\bf H}_{F_1, G_1})$, $({\bf H}^1_{F_2,G_2})$, $({\bf H}^3_{F_2,G_2})$ and $({\bf H}^4_{G_2})$, for $\{u_{\e}, \e\in(0,1)\}\subset\cA_{N}$, there exists a constant $C>0$ independent of $\e, \g$ such that 
 \be
 &&\sup\limits_{t\in[0,T]}\mE|\hat{Y}_{t}^{\e,\g,u_\e}|_\mH^2\leq C(1+|x_0|^{2}+|y_0|^{2}), \label{hatzub}\\
&&\mE\int_0^{T\wedge \tau_1}|Y_{t}^{\e,\g,u_\e}-\hat{Y}_{t}^{\e,\g,u_\e}|_\mH^2\dif t\leq C_{R}\(\d^{1/2}+\frac{\g}{\e}\)\left(1+|x_0|_\mH^2+|y_0|_\mH^2\right).
\label{unztu}
\ee
\el
\begin{proof}
Since the proof of (\ref{hatzub}) is similar to that for (\ref{frozmomeesti}), we only prove (\ref{unztu}).

First of all, applying the It\^{o} formula to $|Y_{t}^{\e,\g,u_\e}-\hat{Y}_{t}^{\e,\g,u_\e}|_\mH^2 e^{l_4 t}$ for any $l_4>0$, one could obtain that
\ce
&&|Y_{t}^{\e,\g,u_\e}-\hat{Y}_{t}^{\e,\g,u_\e}|_\mH^2 e^{l_4 t}\\
&=&l_4\int_{0}^{t}|Y_{s}^{\e,\g,u_\e}-\hat{Y}_{s}^{\e,\g,u_\e}|_\mH^2 e^{l_4 s}\dif s\\
&&+\frac{1}{\g}\int_0^t2e^{l_4 s}{_{\mV^*}}\<AY_{s}^{\e,\g,u_\e}-A\hat{Y}_{s}^{\e,\g,u_\e}, Y_{s}^{\e,\g,u_\e}-\hat{Y}_{s}^{\e,\g,u_\e}\>_\mV\dif s\\
&&+\frac{1}{\g}\int_{0}^{t}2e^{l_4 s}\<F_{2}(X_{s}^{\e,\g,u_\e},Y_{s}^{\e,\g,u_\e})-F_{2}(X_{s(\d)}^{\e,\g,u_\e},\hat{Y}_{s}^{\e,\g,u_\e}), Y_{s}^{\e,\g,u_\e}-\hat{Y}_{s}^{\e,\g,u_\e}\>_\mH\dif s\\
&&+\frac{1}{\sqrt{\g}}\int_{0}^{t}2e^{l_4 s}\<Y_{s}^{\e,\g,u_\e}-\hat{Y}_{s}^{\e,\g,u_\e}, \(G_2(X_s^{\e,\g,u_\e},Y_s^{\e,\g,u_\e})-G_{2}(X_{s(\d)}^{\e,\g,u_\e},\hat{Y}_{s}^{\e,\g,u_\e})\)\dif W^2_{s}\>_\mH\\
&&+\frac{1}{\sqrt{\g \e}}\int_{0}^{t}2e^{l_4 s}\<Y_{s}^{\e,\g,u_\e}-\hat{Y}_{s}^{\e,\g,u_\e}, G_2(X_s^{\e,\g,u_\e},Y_s^{\e,\g,u_\e})\Pi_2u_\e(s)\>_\mH\dif s\\
&&+\int_{0}^{t}\int_0^12e^{l_4 s}(Y_{s}^{\e,\g,u_\e}(\xi)-\hat{Y}_{s}^{\e,\g,u_\e}(\xi))K^{2,\e,\g,u_\e}(\dif s,\dif \xi)\\
&&-\int_{0}^{t}\int_0^12e^{l_4 s}(Y_{s}^{\e,\g,u_\e}(\xi)-\hat{Y}_{s}^{\e,\g,u_\e}(\xi))\hat{K}^{2,\e,\g,u_\e}(\dif s,\dif \xi)\\
&&+\frac{1}{\g}\int_{0}^{t}e^{l_4 s}\|G_2(X_s^{\e,\g,u_\e},Y_s^{\e,\g,u_\e})-G_{2}(X_{s(\d)}^{\e,\g,u_\e},\hat{Y}_{s}^{\e,\g,u_\e})\|^{2}_{\cL_2(\mU_2,\mH)}\dif s.
\de

By the definition of $A$ and the Poincar\'e inequality, we know that
\ce
2{_{\mV^*}}\<AY_{s}^{\e,\g,u_\e}-A\hat{Y}_{s}^{\e,\g,u_\e}, Y_{s}^{\e,\g,u_\e}-\hat{Y}_{s}^{\e,\g,u_\e}\>_\mV\leq -2\l_1|Y_{s}^{\e,\g,u_\e}-\hat{Y}_{s}^{\e,\g,u_\e}|_\mH^2.
\de

Note that
\ce
&&2\<F_{2}(X_{s}^{\e,\g,u_\e},Y_{s}^{\e,\g,u_\e})-F_{2}(X_{s(\d)}^{\e,\g,u_\e},\hat{Y}_{s}^{\e,\g,u_\e}), Y_{s}^{\e,\g,u_\e}-\hat{Y}_{s}^{\e,\g,u_\e}\>_\mH\\
&\leq&2\<F_{2}(X_{s}^{\e,\g,u_\e},Y_{s}^{\e,\g,u_\e})-F_{2}(X_{s}^{\e,\g,u_\e},\hat{Y}_{s}^{\e,\g,u_\e}), Y_{s}^{\e,\g,u_\e}-\hat{Y}_{s}^{\e,\g,u_\e}\>_\mH\\
&&+2\<F_{2}(X_{s}^{\e,\g,u_\e},\hat{Y}_{s}^{\e,\g,u_\e})-F_{2}(X_{s(\d)}^{\e,\g,u_\e},\hat{Y}_{s}^{\e,\g,u_\e}), Y_{s}^{\e,\g,u_\e}-\hat{Y}_{s}^{\e,\g,u_\e}\>_\mH\\
&=:&J_1+J_2.
\de
For $J_1$, by $({\bf H}^1_{F_2, G_2})$, it holds that
\ce
J_1\leq 2L_4|Y_{s}^{\e,\g,u_\e}-\hat{Y}_{s}^{\e,\g,u_\e}|_\mH^2 .
\de
For $J_2$, $({\bf H}^1_{F_2, G_2})$ and the Young inequality imply  that
\ce
J_2&\leq& 2L_3|X_{s}^{\e,\g,u_\e}-X_{s(\d)}^{\e,\g,u_\e}|_\mH |Y_{s}^{\e,\g,u_\e}-\hat{Y}_{s}^{\e,\g,u_\e}|_\mH\\
&\leq&L_4|Y_{s}^{\e,\g,u_\e}-\hat{Y}_{s}^{\e,\g,u_\e}|_\mH^2 +C|X_{s}^{\e,\g,u_\e}-X_{s(\d)}^{\e,\g,u_\e}|_\mH^2.
\de
So,
\ce
&&2\<F_{2}(X_{s}^{\e,\g,u_\e},Y_{s}^{\e,\g,u_\e})-F_{2}(X_{s(\d)}^{\e,\g,u_\e},\hat{Y}_{s}^{\e,\g,u_\e}), Y_{s}^{\e,\g,u_\e}-\hat{Y}_{s}^{\e,\g,u_\e}\>_\mH\\
&\leq&3L_4|Y_{s}^{\e,\g,u_\e}-\hat{Y}_{s}^{\e,\g,u_\e}|_\mH^2 +C|X_{s}^{\e,\g,u_\e}-X_{s(\d)}^{\e,\g,u_\e}|_\mH^2.
\de

By the Young inequality and $({\bf H}^4_{G_2})$, we infer that
\ce
&&\frac{1}{\sqrt{\g \e}}\int_{0}^{t}2e^{l_4 s}\<Y_{s}^{\e,\g,u_\e}-\hat{Y}_{s}^{\e,\g,u_\e}, G_{2}(X_s^{\e,\g,u_\e},Y_s^{\e,\g,u_\e})\Pi_2u_\e(s)\>_\mH\dif s\\
&\leq&\frac{\a}{\g}\int_{0}^{t}e^{l_4 s}|Y_{s}^{\e,\g,u_\e}-\hat{Y}_{s}^{\e,\g,u_\e}|_\mH^2 \dif s+\frac{C}{\e}\int_{0}^{t}e^{l_4 s}\|G_{2}(X_s^{\e,\g,u_\e},Y_s^{\e,\g,u_\e})\|^2_{\cL_2(\mU_2,\mH)}|u_\e(s)|_\sU^2\dif s\\
&\leq&\frac{\a}{\g}\int_{0}^{t}e^{l_4 s}|Y_{s}^{\e,\g,u_\e}-\hat{Y}_{s}^{\e,\g,u_\e}|_\mH^2 \dif s+\frac{C}{\e}L_8^2(1+\sup\limits_{s\in[0,T]}|X_s^{\e,\g,u_\e}|_\mH)^2\int_{0}^{t}e^{l_4 s}|u_\e(s)|_\sU^2\dif s.
\de

Definition \ref{soludefi} yields that
\ce
&&\int_{0}^{t}\int_0^12e^{l_4 s}(Y_{s}^{\e,\g,u_\e}(\xi)-\hat{Y}_{s}^{\e,\g,u_\e}(\xi))K^{2,\e,\g,u_\e}(\dif s,\dif \xi)\leq 0,\\
&&-\int_{0}^{t}\int_0^12e^{l_4 s}(Y_{s}^{\e,\g,u_\e}(\xi)-\hat{Y}_{s}^{\e,\g,u_\e}(\xi))\hat{K}^{2,\e,\g,u_\e}(\dif s,\dif \xi)\leq 0.
\de

$({\bf H}^1_{F_2, G_2})$ implies that
\ce
&&\|G_2(X_s^{\e,\g,u_\e},Y_s^{\e,\g,u_\e})-G_{2}(X_{s(\d)}^{\e,\g,u_\e},\hat{Y}_{s}^{\e,\g,u_\e})\|^{2}_{\cL_2(\mU_2,\mH)}\\
&\leq&2\|G_{2}(X_{s}^{\e,\g,u_\e},Y_{s}^{\e,\g,u_\e})-G_{2}(X_{s}^{\e,\g,u_\e},\hat{Y}_{s}^{\e,\g,u_\e})\|^2_{\cL_2(\mU_2,\mH)}\\
&&+2\|G_{2}(X_{s}^{\e,\g,u_\e},\hat{Y}_{s}^{\e,\g,u_\e})-G_{2}(X_{s(\d)}^{\e,\g,u_\e},\hat{Y}_{s}^{\e,\g,u_\e})\|^2_{\cL_2(\mU_2,\mH)}\\
&\leq&2L_6^2|Y_{s}^{\e,\g,u_\e}-\hat{Y}_{s}^{\e,\g,u_\e}|_\mH^2 +2L_5^2|X_{s}^{\e,\g,u_\e}-X_{s(\d)}^{\e,\g,u_\e}|_\mH^2.
\de

Combining the above deduction, we obtain that
\ce
&&|Y_{t}^{\e,\g,u_\e}-\hat{Y}_{t}^{\e,\g,u_\e}|_\mH^2 e^{l_4 t}\\
&\leq&(l_4-\frac{2\a}{\g}+\frac{\a}{\g})\int_{0}^{t}e^{l_4 s}|Y_{s}^{\e,\g,u_\e}-\hat{Y}_{s}^{\e,\g,u_\e}|_\mH^2 \dif s\\
&&+\frac{1}{\sqrt{\g}}\int_{0}^{t}2e^{l_4 s}\<Y_{s}^{\e,\g,u_\e}-\hat{Y}_{s}^{\e,\g,u_\e}, \(G_2(X_s^{\e,\g,u_\e},Y_s^{\e,\g,u_\e})-G_{2}(X_{s(\d)}^{\e,\g,u_\e},\hat{Y}_{s}^{\e,\g,u_\e})\)\dif W^2_{s}\>_\mH\\
&&+\frac{C}{\g}\int_{0}^{t}e^{l_4 s}|X_{s}^{\e,\g,u_\e}-X_{s(\d)}^{\e,\g,u_\e}|_\mH^2\dif s+\frac{C}{\e}L_8^2(1+\sup\limits_{s\in[0,T]}|X_s^{\e,\g,u_\e}|_\mH)^2\int_{0}^{t}e^{l_4 s}|u_\e(s)|_\sU^2\dif s.
\de
which together with $l_4=\frac{\a}{\g}$ yields that
\ce
&&|Y_{t}^{\e,\g,u_\e}-\hat{Y}_{t}^{\e,\g,u_\e}|_\mH^2 \\
&\leq&\frac{1}{\sqrt{\g}}\int_{0}^{t}2e^{-\frac{\a}{\g}(t-s)}\<Y_{s}^{\e,\g,u_\e}-\hat{Y}_{s}^{\e,\g,u_\e}, \(G_2(X_s^{\e,\g,u_\e},Y_s^{\e,\g,u_\e})-G_{2}(X_{s(\d)}^{\e,\g,u_\e},\hat{Y}_{s}^{\e,\g,u_\e})\)\dif W^2_{s}\>_\mH\\
&&+\frac{C}{\g}\int_{0}^{t}e^{-\frac{\a}{\g}(t-s)}|X_{s}^{\e,\g,u_\e}-X_{s(\d)}^{\e,\g,u_\e}|_\mH^2\dif s+\frac{C}{\e}L_8^2(1+\sup\limits_{s\in[0,T]}|X_s^{\e,\g,u_\e}|_\mH)^2\int_{0}^{t}e^{-\frac{\a}{\g}(t-s)}|u_\e(s)|_\sU^2\dif s.
\de
Integrating from $0$ to $T\wedge \tau_1$ and taking the expectation on two sides, by (\ref{xeguettd}) and (\ref{xeguees}) we conclude that
\ce
&&\mE\int_0^{T\wedge \tau_1}|Y_{t}^{\e,\g,u_\e}-\hat{Y}_{t}^{\e,\g,u_\e}|_\mH^2 \dif t\\
&\leq&\frac{1}{\sqrt{\g}}\mE\int_0^{T\wedge \tau_1}\int_{0}^{t}2e^{-\frac{\a}{\g}(t-s)}\<Y_{s}^{\e,\g,u_\e}-\hat{Y}_{s}^{\e,\g,u_\e}, \\
&&\qquad\qquad\qquad\qquad \(G_2(X_s^{\e,\g,u_\e},Y_s^{\e,\g,u_\e})-G_{2}(X_{s(\d)}^{\e,\g,u_\e},\hat{Y}_{s}^{\e,\g,u_\e})\)\dif W^2_{s}\>_\mH\dif t\\
&&+\frac{C}{\g}\mE\int_0^{T\wedge \tau_1}\int_{0}^{t}e^{-\frac{\a}{\g}(t-s)}|X_{s}^{\e,\g,u_\e}-X_{s(\d)}^{\e,\g,u_\e}|_\mH^2\dif s\dif t\\
&&+\frac{C}{\e}L_8^2\mE(1+\sup\limits_{s\in[0,T]}|X_s^{\e,\g,u_\e}|_\mH)^2\int_0^{T\wedge \tau_1}\int_{0}^{t}e^{-\frac{\a}{\g}(t-s)}|u_\e(s)|_\sU^2\dif s\dif t\\
&\leq&\frac{C}{\a}\mE\int_0^{T\wedge \tau_1}|X_{s}^{\e,\g,u_\e}-X_{s(\d)}^{\e,\g,u_\e}|_\mH^2\dif s+CN\frac{\g}{\e}\mE(1+\sup\limits_{s\in[0,T]}|X_s^{\e,\g,u_\e}|_\mH)^2\\
&\leq&C_{R}\d^{1/2}\left(1+|x_0|_\mH^2+|y_0|_\mH^2\right)+C\frac{\g}{\e}\left(1+|x_0|_\mH^2+|y_0|_\mH^2\right).
\de
The proof is complete.
\end{proof}

Now, we prepare to verify Condition \ref{cond} $(ii)$.

\bp\label{just2}
Assume that $({\bf H}_{F_1, G_1})$, $({\bf H}^1_{F_2, G_2})$, $({\bf H}^2_{G_2})$, $({\bf H}^3_{F_2, G_2})$ and $({\bf H}^4_{G_2})$ hold. Then for each $N<\infty$, any $\{u_\e\}\subset \cA^N$ and any $\eta>0$,
$$
\lim _{\varepsilon \rightarrow 0}\mathbb{P}\left(\sup\limits_{t\in[0,T]}\left|\cG^{\varepsilon}\left(W+\frac{1}{\sqrt{\varepsilon} }\int_0^{\cdot} u_\e(s) \dif s\right)(t)-\cG^0\left(\int_0^{\cdot} u_\e(s)\dif s\right)(t)\right|_\mH>\eta\right)=0.
$$
\ep
\begin{proof}
We divide the proof into three steps. The first step gives the estimate about the difference of $\cG^{\varepsilon}\left(W+\frac{1}{\sqrt{\varepsilon} }\int_0^{\cdot} u_\e(s) \dif s\right)$ and $\cG^0\left(\int_0^{\cdot} u_\e(s)\dif s\right)$ in $C([0,T],\mH)$. In the second step, we deal with the remainder in the first step. The proof is completed in the third step.

{\bf Step 1.} We estimate the difference of $\cG^{\varepsilon}\left(W+\frac{1}{\sqrt{\varepsilon} }\int_0^{\cdot} u_\e(s) \dif s\right)$ and $\cG^0\left(\int_0^{\cdot} u_\e(s)\dif s\right)$ in $C([0,T],\mH)$.

Note that
$$
X^{\e,\g,u_{\e}}=\cG^{\varepsilon}\left(W+\frac{1}{\sqrt{\varepsilon} }\int_0^{\cdot} u_\e(s) \dif s\right), \quad \bar{X}^{u_\e}=\cG^0\left(\int_0^{\cdot} u_\e(s)\dif s\right).
$$
Set $Z^{\e,u_{\e}}_t:=X^{\e,\g,u_{\e}}_{t}-\bar{X}^{u_\e}_{t}$ and by It\^o's formula we get that
\ce
|Z^{\e,u_{\e}}_t|_\mH^2&=&2\int_{0}^t{_{\mV^*}}\<AX^{\e,\g,u_{\e}}_{s}-A\bar{X}^{u_\e}_{s},Z^{\e,u_{\e}}_s\>_\mV\dif s\no\\
&&+2\int_{0}^t{_{\mV^*}}\<B(X^{\e,\g,u_{\e}}_{s},X^{\e,\g,u_{\e}}_{s})-B(\bar{X}^{u_\e}_{s},\bar{X}^{u_\e}_{s}),Z^{\e,u_{\e}}_s\>_\mV\dif s\no\\
&&+2\int_{0}^t\< Z^{\e,u_{\e}}_s, F_1(X^{\e,\g,u_{\e}}_{s},Y^{\e,\g,u_{\e}}_{s})-\bar{F}_1(\bar{X}^{u_\e}_{s})\>_\mH\dif s   \no\\
&&+2\int_{0}^t\<Z^{\e,u_{\e}}_s, G_1(X^{\e,\g,u_{\e}}_{s})\Pi_1u_{\e}(s)-G_1(\bar{X}^{u_\e}_{s})\Pi_1u_\e(s)\>_\mH \dif s  \no\\
&&+2\int_{0}^t\int_0^1Z^{\e,u_{\e}}_s(\xi)K^{1,\e,\g,u_{\e}}(\dif s,\dif \xi)-2\int_{0}^{t}\int_0^1Z^{\e,u_{\e}}_s(\xi)\bar{K}^{1,u_\e}(\dif s,\dif\xi)\no\\
&&+2\sqrt{\e} \int_{0}^t\<Z^{\e,u_{\e}}_s,  G_1(X^{\e,\g,u_{\e}}_{s})\dif W^1_s\>_\mH +\e\int_{0}^{t}\|G_1(X^{\e,\g,u_{\e}}_{s})\|_{\cL_2(\mU_1,\mH)}^{2} \dif s.
\de

By the definition of $A$, it holds that
\ce
2{_{\mV^*}}\<AX^{\e,\g,u_{\e}}_{s}-A\bar{X}^{u_\e}_{s},Z^{\e,u_{\e}}_s\>_\mV=-2\|Z^{\e,u_{\e}}_s\|_\mV^2.
\de
And Lemma \ref{Bbprop3} implies that
\ce
&&2{_{\mV^*}}\<B(X^{\e,\g,u_{\e}}_{s},X^{\e,\g,u_{\e}}_{s})-B(\bar{X}^{u_\e}_{s},\bar{X}^{u_\e}_{s}),Z^{\e,u_{\e}}_s\>_\mV\\
&\leq& 4|Z^{\e,u_{\e}}_s|_\mH(\|X^{\e,\g,u_{\e}}_{s}\|_\mV+\|\bar{X}^{u_\e}_{s}\|_\mV)\|Z^{\e,u_{\e}}_s\|_\mV\\
&\leq&8|Z^{\e,u_{\e}}_s|_\mH^2(\|X^{\e,\g,u_{\e}}_{s}\|^2_\mV+\|\bar{X}^{u_\e}_{s}\|^2_\mV)+\|Z^{\e,u_{\e}}_s\|_\mV^2.
\de
Next, note that
\ce
&&2\< Z^{\e,u_{\e}}_s, F_1(X^{\e,\g,u_{\e}}_{s},Y^{\e,\g,u_{\e}}_{s})-\bar{F}_1(\bar{X}^{u_\e}_{s})\>_\mH\\
&=&\bigg[2\<Z^{\e,u_{\e}}_s, F_1(X^{\e,\g,u_{\e}}_{s},Y^{\e,\g,u_{\e}}_{s})-F_1(X^{\e,\g,u_{\e}}_{s(\d)},\hat{Y}^{\e,\g,u_{\e}}_{s}) \>_\mH\\
&&+2\<Z^{\e,u_{\e}}_s,-\bar{F}_1(X^{\e,\g,u_{\e}}_{s})+\bar{F}_1(X^{\e,\g,u_{\e}}_{s(\d)})\>_\mH\\
&&+2\<Z^{\e,u_{\e}}_s,\bar{F}_1(X^{\e,\g,u_{\e}}_{s})-\bar{F}_1(\bar{X}^{u_\e}_{s}) \>_\mH\bigg]\\
&&+2\<Z^{\e,u_{\e}}_s-Z^{\e,u_{\e}}_{s(\d)},F_1(X^{\e,\g,u_{\e}}_{s(\d)},\hat{Y}^{\e,\g,u_{\e}}_{s})-\bar{F}_1(X^{\e,\g,u_{\e}}_{s(\d)})\>_\mH\\
&&+2\<Z^{\e,u_{\e}}_{s(\d)},F_1(X^{\e,\g,u_{\e}}_{s(\d)},\hat{Y}^{\e,\g,u_\e}_{s})-\bar{F}_1(X^{\e,\g,u_{\e}}_{s(\d)})\>_\mH\\
&=:&\cK_1(s)+\cK_2(s)+\cK_3(s).
\de
So, by the Lipschitz continuity of $F_1, \bar{F}_1$ and the H\"older inequality, we get that
\ce
\int_{0}^t\cK_1(s)\dif s&\leq& C\int_{0}^t|Z^{\e,u_{\e}}_s|_\mH^2\dif s+C\int_{0}^t|X^{\e,\g,u_{\e}}_{s}-X^{\e,\g,u_{\e}}_{s(\d)}|_\mH^2\dif s\\
&&+C\int_{0}^t|Y^{\e,\g,u_{\e}}_{s}-\hat{Y}^{\e,\g,u_{\e}}_{s}|_\mH^2\dif s,
\de
and
\ce
\int_{0}^t\cK_2(s)\dif s\leq C\int_{0}^t(|X^{\e,\g,u_{\e}}_s-X^{\e,\g,u_{\e}}_{s(\d)}|_\mH+|\bar{X}^{u_\e}_s-\bar{X}^{u_\e}_{s(\d)}|_\mH)(1+|X^{\e,\g,u_{\e}}_{s(\d)}|_\mH+|\hat{Y}^{\e,\g,u_{\e}}_{s}|_\mH)\dif s.
\de

Besides, it is easy to see that
\ce
&&2\int_{0}^t\<Z^{\e,u_{\e}}_s, G_1(X^{\e,\g,u_{\e}}_{s})\Pi_1u_{\e}(s)-G_1(\bar{X}^{u_\e}_{s})\Pi_1u_\e(s)\>_\mH\dif s\\
&\leq&2\int_{0}^t\left|Z^{\e,u_{\e}}_s\right|_\mH\cdot\left|(G_1(X^{\e,\g,u_{\e}}_{s})-G_1(\bar{X}^{u_\e}_{s}))\Pi_1u_\e(s)\right|_\mH\dif s\\
&\leq&2L_2\int_{0}^t\left|Z^{\e,u_{\e}}_s\right|_\mH^2|u_\e(s)|_\sU\dif s.
\de
And Definition \ref{soludefi} yields that
\ce
2\int_{0}^t\int_0^1Z^{\e,u_{\e}}_s(\xi)K^{1,\e,\g,u_{\e}}(\dif s,\dif \xi)-2\int_{0}^t\int_0^1Z^{\e,u_{\e}}_s(\xi)\bar{K}^{1,u_\e}(\dif s,\dif\xi)\leq 0.
\de

Combining the above deduction, we obtain that
\ce
&&|Z^{\e,u_{\e}}_t|^{2}_{\mH}+\int_{0}^t\|Z^{\e,u_{\e}}_s\|_\mV^2\dif s\\
&\leq& C\int_{0}^t|Z^{\e,u_{\e}}_s|_\mH^2(1+|u_\e(s)|_\sU+\|X^{\e,\g,u_{\e}}_{s}\|^2_\mV+\|\bar{X}^{u_\e}_{s}\|^2_\mV)\dif s\\
&&+C\int_{0}^t|X^{\e,\g,u_{\e}}_{s}-X^{\e,\g,u_{\e}}_{s(\d)}|_\mH^2\dif s+C\int_{0}^t|Y^{\e,\g,u_{\e}}_{s}-\hat{Y}^{\e,\g,u_{\e}}_{s}|_\mH^2\dif s\\
&&+\int_{0}^t(|X^{\e,\g,u_{\e}}_s-X^{\e,\g,u_{\e}}_{s(\d)}|_\mH+|\bar{X}^{u_\e}_s-\bar{X}^{u_\e}_{s(\d)}|_\mH)(1+|X^{\e,\g,u_{\e}}_{s(\d)}|_\mH+|\hat{Y}^{\e,\g,u_{\e}}_{s}|_\mH)\dif s\\
&&+2\left|\int_{0}^t\<Z^{\e,u_{\e}}_{s(\d)},F_1(X^{\e,\g,u_{\e}}_{s(\d)},\hat{Y}^{\e,\g,u_\e}_{s})-\bar{F}_1(X^{\e,\g,u_{\e}}_{s(\d)})\>_\mH\dif s\right|\\
&&+2\sqrt{\e} \left|\int_{0}^t\<Z^{\e,u_{\e}}_s,  G_1(X^{\e,\g,u_{\e}}_{s})\dif W^1_s\>_\mH\right|+\e\int_{0}^{t}\|G_1(X^{\e,\g,u_{\e}}_{s})\|_{\cL_2(\mU_1,\mH)}^{2} \dif s.
\de
Besides, we define the stopping time 
\be
\tau_2:=\inf\left\{t>0, |X_t^{\e,\g,u_\e}|_\mH+\int_0^t\(\|X^{\e,\g,u_{\e}}_{s}\|_\mV^2+\|\bar{X}^{u_\e}_{s}\|^2_\mV\)\dif s>R\right\},
\label{tau2defi}
\ee
and $\tau_2\leq \tau_1$. So, by the Gronwall inequality it holds that
\be
&&\mE\sup\limits_{t\in[0,T\wedge \tau_2]}|Z^{\e,u_{\e}}_t|_\mH^2+\mE\int_{0}^{T\wedge \tau_2}\|Z^{\e,u_{\e}}_s\|_\mV^2\dif s\no\\
&\leq& C_R\mE\int_{0}^{T\wedge \tau_2}|X^{\e,\g,u_{\e}}_{s}-X^{\e,\g,u_{\e}}_{s(\d)}|_\mH^2\dif s+C_R\mE\int_{0}^{T\wedge \tau_2}|Y^{\e,\g,u_{\e}}_{s}-\hat{Y}^{\e,\g,u_{\e}}_{s}|_\mH^2\dif s\no\\
&&+C_R\mE\int_0^{T\wedge \tau_2}(|X^{\e,\g,u_{\e}}_s-X^{\e,\g,u_{\e}}_{s(\d)}|_\mH+|\bar{X}^{u_\e}_s-\bar{X}^{u_\e}_{s(\d)}|_\mH)(1+|X^{\e,\g,u_{\e}}_{s(\d)}|_\mH+|\hat{Y}^{\e,\g,u_{\e}}_{s}|_\mH)\dif s\no\\
&&+2C_R\mE\sup\limits_{t\in[0,T\wedge \tau_2]}\left|\int_{0}^t\<Z^{\e,u_{\e}}_{s(\d)},F_1(X^{\e,\g,u_{\e}}_{s(\d)},\hat{Y}^{\e,\g,u_\e}_{s})-\bar{F}_1(X^{\e,\g,u_{\e}}_{s(\d)})\>_\mH\dif s\right|\no\\
&&+2\sqrt{\e}C_R\mE\sup\limits_{t\in[0,T\wedge \tau_2]}\left|\int_{0}^t\<Z^{\e,u_{\e}}_s,  G_1(X^{\e,\g,u_{\e}}_{s})\dif W^1_s\>_\mH\right|\no\\
&&+\e C_R\mE\int_{0}^{{T\wedge \tau_2}}\|G_1(X^{\e,\g,u_{\e}}_{s})\|_{\cL_2(\mU_1,\mH)}^{2} \dif s\no\\
&=:&J_1+J_2+J_3+J_4+J_5+J_6.
\label{j123456}
\ee

In the following, (\ref{xeguettd}) and (\ref{unztu}) imply that
\be
J_1+J_2\leq C_{R}\(\frac{\g}{\e}+\d^{1/2}\)\left(1+|x_0|_\mH^2+|y_0|_\mH^2\right).
\label{j12}
\ee
Then by (\ref{xeguettd}), (\ref{barxueincr}), (\ref{xeguees}), (\ref{hatzub}) and the H\"older inequality, it holds that
\be
J_3&\leq& C_R\left(\mE\int_0^{T\wedge \tau_2}(|X^{\e,\g,u_{\e}}_s-X^{\e,\g,u_{\e}}_{s(\d)}|_\mH^2+|\bar{X}^{u_\e}_s-\bar{X}^{u_\e}_{s(\d)}|_\mH^2)\dif s\right)^{1/2}\no\\
&&\times\left(\mE\int_0^{T\wedge \tau_2}(1+|X^{\e,\g,u_{\e}}_{s(\d)}|_\mH^2+|\hat{Y}^{\e,\g,u_{\e}}_{s}|_\mH^2)\dif s\right)^{1/2}\no\\
&\leq&C_{R}\d^{1/4}\left(1+|x_0|_\mH^2+|y_0|_\mH^2\right).
\label{j3}
\ee

By the deduction in {\bf Step 2}, we infer that
\be
J_4\leq C_R\((\frac{\g}{\d})^{1/2}+\d^{1/2}\)\left(1+|x_0|_\mH^2+|y_0|_\mH^2\right).
\label{j4}
\ee

For $J_5$, from the Burkholder-Davis-Gundy inequality, the Young inequality and the linear growth of $G_1$, it follows that
\be
J_5&\leq& 2\sqrt{\e}C\mE\left(\int_{0}^{T\wedge \tau_2}|Z^{\e,u_{\e}}_s|_\mH^2\|G_1(X^{\e,\g,u_{\e}}_{s}) \|_{\cL_2(\mU_1,\mH)}^2 \dif s\right)^{1/2}\no\\
&\leq&\frac{1}{2}\mE\sup\limits_{t\in[0,T\wedge \tau_2]}|Z^{\e,u_{\e}}_t|_\mH^{2}+C\e\mE\int_{0}^{T\wedge \tau_2}(1+|X^{\e,\g,u_{\e}}_{s}|_\mH^2)\dif s\no\\
&\leq&\frac{1}{2}\mE\sup\limits_{t\in[0,T\wedge \tau_2]}|Z^{\e,u_{\e}}_t|_\mH^{2}+CT\e\left(1+|x_0|_\mH^2+|y_0|_\mH^2\right).
\label{j5}
\ee
For $J_6$, by the linear growth of $G_1$, we know
\be
J_6\leq \e C_R\mE\int_{0}^{{T\wedge \tau_2}}(1+|X^{\e,\g,u_{\e}}_{s}|_\mH^2)\dif s\leq \e C_R\left(1+|x_0|_\mH^2+|y_0|_\mH^2\right).
\label{j6}
\ee

Combining (\ref{j12})-(\ref{j6}) with (\ref{j123456}), we can get
\be
&&\mE\sup\limits_{t\in[0,T\wedge \tau_2]}|Z^{\e,u_{\e}}_t|_\mH^{2}+\mE\int_{0}^{T\wedge \tau_2}\|Z^{\e,u_{\e}}_s\|_\mV^2\dif s\no\\
&\leq&C_{R}\left(\frac{\g}{\e}+\d^{1/2}+\d^{1/4}+\e+(\frac{\g}{\d})^{1/2}\right)\left(1+|x_0|_\mH^2+|y_0|_\mH^2\right).
\label{supintes}
\ee

{\bf Step 2.} We prove (\ref{j4}).

For $J_4$, it holds that
\ce
J_4&\leq&2C_R\mE\sup\limits_{t\in[0,T\wedge \tau_2]}\left|\int_{0}^{[\frac{t}{\d}]\d}\<Z^{\e,u_{\e}}_{s(\d)},F_1(X^{\e,\g,u_{\e}}_{s(\d)},\hat{Y}^{\e,\g,u_\e}_{s})-\bar{F}_1(X^{\e,\g,u_{\e}}_{s(\d)})\>_\mH\dif s\right|\no\\
&&+2C_R\mE\sup\limits_{t\in[0,T\wedge \tau_2]}\left|\int_{[\frac{t}{\d}]\d}^t\<Z^{\e,u_{\e}}_{s(\d)},F_1(X^{\e,\g,u_{\e}}_{s(\d)},\hat{Y}^{\e,\g,u_\e}_{s})-\bar{F}_1(X^{\e,\g,u_{\e}}_{s(\d)})\>_\mH\dif s\right|\no\\
 &=:&J_{41}+J_{42}.
\de

Next, we are devoted to estimating $J_{41}$. Note that
\be
J_{41}
&=&2C_R\mE\sup\limits_{t\in[0,T\wedge \tau_2]}
\Big|\sum\limits_{k=0}^{[\frac{t}{\d}]-1}\int_{k\d}^{(k+1)\d}\<Z^{\e,u_{\e}}_{s(\d)},F_1(X^{\e,\g,u_{\e}}_{s(\d)},\hat{Y}^{\e,\g,u_\e}_{s})-\bar{F}_1(X^{\e,\g,u_{\e}}_{s(\d)})\>_\mH\dif s\Big|\no\\
&\leq&2C_R\mE\sup\limits_{t\in[0,T\wedge \tau_2]}\sum\limits_{k=0}^{[\frac{t}{\d}]-1}\left|\int_{k\d}^{(k+1)\d}\<Z^{\e,u_{\e}}_{k\d},F_1(X^{\e,\g,u_{\e}}_{k\d},\hat{Y}^{\e,\g,u_\e}_{s})-\bar{F}_1(X^{\e,\g,u_{\e}}_{k\d})\>_\mH\dif s\right|\no\\
&\leq&2C_R\sum\limits_{k=0}^{[\frac{T}{\d}]-1}\mE\Big|\int_{k\d}^{(k+1)\d}\<Z^{\e,u_{\e}}_{k\d},F_1(X^{\e,\g,u_{\e}}_{k\d},\hat{Y}^{\e,\g,u_\e}_{s})-\bar{F}_1(X^{\e,\g,u_{\e}}_{k\d})\>_\mH\dif s\Big|\no\\
&\leq&2C_R[\frac{T}{\d}]\sup_{0\leq k\leq [\frac{T}{\d}]-1}\mE\Big|\int_{k\d}^{(k+1)\d}\<Z^{\e,u_{\e}}_{k\d},F_1(X^{\e,\g,u_{\e}}_{k\d},\hat{Y}^{\e,\g,u_\e}_{s})-\bar{F}_1(X^{\e,\g,u_{\e}}_{k\d})\>_\mH\dif s\Big|\no\\
&\leq&2C_R\g(\frac{T}{\d})\sup_{0\leq k\leq [\frac{T}{\d}]-1}\mE\Big|\<Z^{\e,u_{\e}}_{k\d},\int_{0}^{\d/\g}(F_1(X^{\e,\g,u_{\e}}_{k\d},\hat{Y}^{\e,\g,u_\e}_{s\g+k\d})-\bar{F}_1(X^{\e,\g,u_{\e}}_{k\d}))\dif s\>_\mH\Big|\no\\
&\leq&2C_R\g(\frac{T}{\d})\sup_{0\leq k\leq [\frac{T}{\d}]-1}\mE|Z^{\e,u_{\e}}_{k\d}|_\mH\left|\int_{0}^{\d/\g}(F_1(X^{\e,\g,u_{\e}}_{k\d},\hat{Y}^{\e,\g,u_\e}_{s\g+k\d})-\bar{F}_1(X^{\e,\g,u_{\e}}_{k\d}))\dif s\right|_\mH\no\\
 &\leq&2C_R\g(\frac{T}{\d})\sup_{0\leq k\leq [\frac{T}{\d}]-1}(\mE|Z^{\e,u_{\e}}_{k\d}|_\mH^2)^{1/2}\no\\
 &&\times\left(\mE\left|\int_{0}^{\d/\g}(F_1(X^{\e,\g,u_{\e}}_{k\d},\hat{Y}^{\e,\g,u_\e}_{s\g+k\d})-\bar{F}_1(X^{\e,\g,u_{\e}}_{k\d}))\dif s\right|_\mH^2\right)^{1/2}\no\\
  &\leq&2C_R\g(\frac{T}{\d})\left(1+|x_0|_\mH^2+|y_0|_\mH^2\right)^{1/2}\sup_{0\leq k\leq [\frac{T}{\d}]-1}\left(2\int_{0}^{\d/\g}\int_{r}^{\d/\g}\Phi(s,r)\dif s\dif r\right)^{1/2},
\label{b41c}
\ee
where for $0\leq r\leq s\leq\d/\g$
$$
\Phi(s,r):=\mE\<F_1(X^{\e,\g,u_{\e}}_{k\d},\hat{Y}^{\e,\g,u_\e}_{s\g+k\d})-\bar{F}_1(X^{\e,\g,u_{\e}}_{k\d}),F_1(X^{\e,\g,u_{\e}}_{k\d},\hat{Y}^{\e,\g,u_\e}_{r\g+k\d})-\bar{F}_1(X^{\e,\g,u_{\e}}_{k\d})\>_\mH.
$$

In the following, for any $w>0$ and random variables $X, Y\in L^2(\Omega,\mathscr{F}_w,\mP; \mV)$, we construct the following SPDE with reflection
\ce\left\{\begin{array}{l}
\dif \check{Y}_t^{\g, X,Y}=\frac{1}{\g}[A\check{Y}_t^{\g,X,Y}+F_2(X, \check{Y}_t^{\g,X,Y})]\dif t+\frac{1}{\sqrt{\g}}G_2(X, \check{Y}_t^{\g,X,Y})\dif W^2_t+\dif \check{K}_t^{2,\g,X,Y},\quad t\geq w,\\
\check{Y}_{t}^{\g,X,Y}(\xi)\geq 0,\quad \xi\in[0,1],\\
\check{Y}_{t}^{\g,X,Y}(0)=\check{Y}_{t}^{\g,X,Y}(1)=0,\\
\check{Y}_w^{\g,X,Y}=Y.
\end{array}
\right.
\de
Under $({\bf H}^1_{F_2, G_2})$, the above SPDE with reflection has a unique solution $(\check{Y}^{\g,X,Y}, \check{K}^{2,\g,X,Y})$. Moreover, it holds that for $t \in[k \delta,(k+1) \delta]$
\ce
&&\hat{Y}_t^{\e,\g,u_\e}=\check{Y}_t^{\g,X_{k \delta}^{\e,\g,u_\e}, \hat{Y}_{k \delta}^{\e,\g,u_\e}},\\
&&\hat{K}_t^{2,\e,\g,u_\e}=\check{K}_t^{2,\g,X_{k \delta}^{\e,\g,u_\e}, \hat{Y}_{k \delta}^{\e,\g,u_\e}},
\de
and
\ce
\Phi(s,r)&=&\mE\<F_1(X^{\e,\g,u_{\e}}_{k\d},\check{Y}^{\g,X_{k \delta}^{\e,\g,u_\e}, \hat{Y}_{k \delta}^{\e,\g,u_\e}}_{\g s+k\d})-\bar{F}_1(X^{\e,\g,u_{\e}}_{k\d}),\\
&&\qquad F_1(X^{\e,\g,u_{\e}}_{k\d},\check{Y}^{\g,X_{k \delta}^{\e,\g,u_\e}, \hat{Y}_{k \delta}^{\e,\g,u_\e}}_{\g r+k\d})-\bar{F}_1(X^{\e,\g,u_{\e}}_{k\d})\>_\mH.
\de
Since $X_{k \delta}^{\e,\g,u_\e}, \hat{Y}_{k \delta}^{\e,\g,u_\e}$ are $\sF_{k\d}$-measurable, and for any $x,y\in\mV$, $\check{Y}_t^{\g,x,y}$ is independent of $\sF_{k\d}$, we obtain that
\ce
\Phi(s,r)
 &=&\mE\Bigg[\mE\Bigg[\<F_{1}(X_{k\d}^{\e,\g,u_\e},\check{Y}^{\g, X_{k \delta}^{\e,\g,u_\e}, \hat{Y}_{k \delta}^{\e,\g,u_\e}}_{\g s+k\d})
 -\bar{F}_{1}(X_{k\d}^{\e,\g,u_\e}),\\
 &&\qquad\qquad F_{1}(X_{k\d}^{\e,\g,u_\e},\check{Y}^{\g,X_{k \delta}^{\e,\g,u_\e}, \hat{Y}_{k \delta}^{\e,\g,u_\e}}_{\g r+k\d})
 -\bar{F}_{1}(X_{k\d}^{\e,\g,u_\e})\>_\mH\Bigg{|}\sF_{k\d}\Bigg]\Bigg]\\
&=&\mE\Bigg[\mE\Bigg[\<F_{1}(x,\check{Y}^{\g,x, y}_{\g s+k\d})
 -\bar{F}_{1}(x),F_{1}(x,\check{Y}^{\g,x, y}_{\g r+k\d})
 -\bar{F}_{1}(x)\>_\mH\Bigg]\Bigg{|}_{(x,y)=(X_{k\d}^{\e,\g,u_\e},\hat{Y}_{k \delta}^{\e,\g,u_\e})}\Bigg].
\de

Now, we investigate $\check{Y}^{\g,x, y}_{\g s+k\d}$. On the one hand, it holds that
\ce
\check{Y}^{\g, x, y}_{\g s+k\d}&=&y+\frac{1}{\g} \int_{k\d}^{\g s+k\d}[A\check{Y}^{\g,x, y}_{r}+F_2(x, \check{Y}^{\g,x, y}_{r})]\dif r+\frac{1}{\sqrt{\g}} \int_{k\d}^{\g s+k\d} G_2(x, \check{Y}^{\g,x, y}_{r})\dif W^2_r\\
&&+\check{K}_{\g s+k\d}^{2,\g,x,y}-\check{K}_{k\d}^{2,\g,x,y}\\
&=&y+\frac{1}{\g} \int_{0}^{\g s}\left[A\check{Y}^{\g,x, y}_{\iota+k\d}+F_2(x, \check{Y}^{\g,x, y}_{\iota+k\d})\right]\dif \iota+\frac{1}{\sqrt{\g}} \int_{0}^{\g s} G_2(x, \check{Y}^{\g,x, y}_{\iota+k\d})\dif \hat{W}^2_\iota\\
&&+\check{K}_{\g s+k\d}^{2,\g,x,y}-\check{K}_{k\d}^{2,\g, x,y}\\
&=&y+\int_{0}^{s}\left[A\check{Y}^{\g,x, y}_{\g v+k\d}+F_2(x,\check{Y}^{\g,x, y}_{\g v+k\d})\right]\dif v+\int_{0}^{s} G_2(x, \check{Y}^{\g,x, y}_{\g v+k\d})\dif \check{\hat{W}}^2_v+\check{\check{K}}_{s}^{2,\g,x,y},
\de
where $\hat{W}^2_\cdot:=W^2_{\cdot+k\d}-W^2_{k\d}$ and $\check{\hat{W}}^2_\cdot:=\frac{1}{\sqrt{\g}}\hat{W}^2_{\g \cdot}$ are two $\mU_2$-valued cylindrical Wiener processes, and $\check{\check{K}}_{s}^{2,\g,x,y}:=\check{K}_{\g s+k\d}^{2,\g,x,y}-\check{K}_{k\d}^{2,\g,x,y}$. On the other hand, the frozen equation (\ref{frozequa}) is written as
\ce
Y_{s}^{x, y}=y+\int_0^s \left[AY_{r}^{x, y}+F_{2}(x,Y_{r}^{x, y})\right]\dif r+\int_0^s G_2(x,Y_{r}^{x, y})\dif \tilde W^2_{r}+K^{2,x, y}_s.
\de
Thus, for $s\in[0,\d/\g]$, $\check{Y}^{\g,x, y}_{\g s+k\d}$ and $Y_{s}^{x, y}$ have the same distribution, which implies that
\ce
&&\mE\Bigg[\<F_{1}(x,\check{Y}^{\g,x, y}_{\g s+k\d})-\bar{F}_{1}(x),F_{1}(x,\check{Y}^{\g,x, y}_{\g r+k\d})-\bar{F}_{1}(x)\>_\mH\Bigg]\\
&=&\mE\<F_{1}(x,Y_{s}^{x, y})-\bar{F}_{1}(x),F_{1}(x,Y_{r}^{x, y})-\bar{F}_{1}(x)\>_\mH\\
&=&\mE\left[\mE\left[\Bigg<F_{1}(x,Y_{s}^{x, y})-\bar{F}_{1}(x), F_{1}(x,Y_{r}^{x, y})-\bar{F}_{1}(x)\Bigg>_\mH\Bigg{|}\sF_r^{\tilde W^2}\right]\right]\\
&=&\mE\left[\Bigg<\mE\left[F_{1}(x,Y_{s}^{x, y})\Bigg{|}\sF_r^{\tilde W^2}\right]-\bar{F}_{1}(x),F_{1}(x,Y_{r}^{x, y})-\bar{F}_{1}(x)\Bigg>_\mH\right]\\
&\leq&\left(\mE\left|\mE\left[F_{1}(x,Y_{s-r}^{x,\hat{y}})\right]\Big{|}_{\hat{y}=Y_{r}^{x,y}}-\bar{F}_{1}(x)\right|_\mH^2\right)^{1/2}\left(\mE|F_{1}(x,Y_{r}^{x, y})-\bar{F}_{1}(x)|_\mH^2\right)^{1/2},
\de
where $(\mathscr{F}_{t}^{\tilde W^2})_{t\in[0,\infty)}$ is the natural filtration generated by $(\tilde W^2_{t})$. 

Besides, by (\ref{meu2}) and (\ref{frozmomeesti}) it holds that
\ce
\left(\mE\left|\mE\left[F_{1}(x,Y_{s-r}^{x,\hat{y}})\right]\Big{|}_{\hat{y}=Y_{r}^{x,y}}-\bar{F}_{1}(x)\right|_\mH^2\right)^{1/2}
&\leq&\(Ce^{-2\a(s-r)}(1+|x|_\mH^2+\mE|Y_{r}^{x,y}|_\mH^2)\)^{\frac{1}{2}}\\
&\leq&\(Ce^{-2\a(s-r)}(1+|x|_\mH^2+|y|_\mH^2)\)^{\frac{1}{2}}\\
&\leq& Ce^{-\a(s-r)}(1+|x|_\mH+|y|_\mH),
\de
and
\ce
\left(\mE|F_{1}(x,Y_{r}^{x, y})-\bar{F}_{1}(x)|_\mH^2\right)^{1/2}
&=&\left(\mE|F_{1}(x,Y_{r}^{x,y})-\int_{\mH}F_{1}(x,z)\mu^{x}(\dif z)|_\mH^2\right)^{1/2}\\
&\leq&\left(\mE\int_{\mH}|F_{1}(x,Y_{r}^{x,y})-F_{1}(x,z)|_\mH^2\mu^{x}(\dif z)\right)^{1/2}\\
&\leq&L_1\left(\int_{\mH}\mE|Y_{r}^{x,y}-z|_\mH^2\mu^{x}(\dif z)\right)^{1/2}\\
&\leq&C\left(|y|_\mH^2 e^{-2\a r}+C(1+|x|_\mH^2)\right)^{1/2}\\
&\leq&C(1+|x|_\mH+|y|_\mH).
\de

The above deduction with (\ref{xeguees}) and (\ref{hatzub}) yields that
\ce
\Phi(s,r)\leq Ce^{-\a(s-r)}\left(1+|x_0|_\mH^2+|y_0|_\mH^2\right).
\de
Inserting the above inequality in (\ref{b41c}), we get that
\be
J_{41}&\leq& 2C_R\g(\frac{T}{\d})\left(1+|x_0|_\mH^2+|y_0|_\mH^2\right)\sup_{0\leq k\leq [\frac{T}{\d}]-1}\left(\int_{0}^{\frac{\d}{\g}}\int_{r}^{\frac{\d}{\g}}Ce^{-\a(s-r)}\dif s\dif r\right)^{1/2}\no\\
&\leq& C_R(\frac{\g}{\d})^{1/2}\left(1+|x_0|_\mH^2+|y_0|_\mH^2\right).
\label{b4de1}
\ee

Next, we deal with $J_{42}$. By (\ref{f1g1linegrow}), (\ref{xeguees}), (\ref{barxue}), (\ref{hatzub}) and the H\"older inequality, one could get that
\be
J_{42}&\leq&2C_R\mE\sup\limits_{t\in[0,T\wedge \tau_2]}\int_{[\frac{t}{\d}]\d}^t|Z^{\e,u_{\e}}_{s(\d)}|_\mH|F_1(X^{\e,\g,u_{\e}}_{s(\d)},\hat{Y}^{\e,\g,u_\e}_{s})-\bar{F}_1(X^{\e,\g,u_{\e}}_{s(\d)})|_\mH\dif s\no\\
 &\leq&2C_R\d^{1/2}\left(\mE\sup_{t\in[0,T]}\int_{[\frac{t}{\d}]\d}^{t}|F_{1}(X_{s(\d)}^{\e,\g,u_\e},\hat{Y}_{s}^{\e,\g,u_\e})-\bar{F}_{1}(X_{s(\d)}^{\e,\g,u_\e})|_\mH^2\dif s\right)^{1/2}\no\\
 &&\times\left(\mE\sup\limits_{s\in[0,T]}|Z^{\e,u_{\e}}_s|_\mH^2\right)^{1/2}\no\\ 
 &\leq&2C_R\d^{1/2}\left(\mE\int_{0}^{T}|F_{1}(X_{s(\d)}^{\e,\g,u_\e},\hat{Y}_{s}^{\e,\g,u_\e})-\bar{F}_{1}(X_{s(\d)}^{\e,\g,u_\e})|_\mH^2\dif s\right)^{1/2}\no\\
 &&\times\left(\mE\sup\limits_{s\in[0,T]}|Z^{\e,u_{\e}}_s|_\mH^2\right)^{1/2}\no\\ 
 &\leq&C_R\d^{1/2}\left(\int_{0}^{T}(1+\mE|X_{s(\d)}^{\e,\g,u_\e}|_\mH^2+\mE|\hat{Y}_{s}^{\e,\g,u_\e}|_\mH^2)\dif s\right)^{1/2}\no\\
 &&\times\left(\mE\sup_{s\in[0,T]}|X^{\e,\g,u_{\e}}_{s}|_\mH^2+\mE\sup_{s\in[0,T]}|\bar{X}^{u_\e}_{s}|_\mH^2\right)^{1/2}\no\\
&\leq& C_R\d^{1/2}\left(1+|x_0|_\mH^2+|y_0|_\mH^2\right).
\label{j42}
\ee

Finally, combining (\ref{b4de1}) and (\ref{j42}) we conclude (\ref{j4}).

{\bf Step 3.} We prove that for each $N<\infty$, any $\{u_\e\}\subset \cA^N$ and any $\eta>0$,
\ce
\lim _{\varepsilon \rightarrow 0}\mP\left(\sup\limits_{t\in[0,T]}\left|X^{\e,\g,u_{\e}}_t-\bar{X}^{u_\e}_{t}\right|_\mH>\eta\right)=0.
\de

First of all, it is easy to see that
\ce
&&\mP\left(\sup\limits_{t\in[0,T]}\left|X^{\e,\g,u_{\e}}_t-\bar{X}^{u_\e}_{t}\right|_\mH>\eta\right)\\
&=&\mP\left(\sup\limits_{t\in[0,T]}\left|X^{\e,\g,u_{\e}}_t-\bar{X}^{u_\e}_{t}\right|_\mH>\eta, T>\tau_2\right)+\mP\left(\sup\limits_{t\in[0,T]}\left|X^{\e,\g,u_{\e}}_t-\bar{X}^{u_\e}_{t}\right|_\mH>\eta, T\leq\tau_2\right)\\
&\leq&\mP\left(|X_T^{\e,\g,u_\e}|_\mH+\int_0^T\(\|X^{\e,\g,u_{\e}}_{s}\|_\mV^2+\|\bar{X}^{u_\e}_{s}\|^2_\mV\)\dif s>R\right)\\
&&+\mP\left(\sup\limits_{t\in[0,T\wedge\tau_2]}\left|X^{\e,\g,u_{\e}}_t-\bar{X}^{u_\e}_{t}\right|_\mH>\eta\right)\\
&\leq&\frac{1}{R}\left[\mE|X_T^{\e,\g,u_\e}|+\mE\int_0^T\(\|X^{\e,\g,u_{\e}}_{s}\|_\mV^2+\|\bar{X}^{u_\e}_{s}\|^2_\mV\)\dif s\right]+\frac{1}{\eta^2}\mE\sup\limits_{t\in[0,T\wedge\tau_2]}\left|X^{\e,\g,u_{\e}}_t-\bar{X}^{u_\e}_{t}\right|_\mH^2\\
&=:&\cM_1+\cM_2.
\de
For $\cM_1$, (\ref{xeguees}) and (\ref{barxue}) imply that
\ce
\cM_1\leq \frac{C}{R}\left(1+|x_0|_\mH^2+|y_0|_\mH^2\right).
\de
For $\cM_2$, by (\ref{supintes}), it holds that
\ce
\cM_2\leq\frac{C_{R}}{\eta^2}\left(\frac{\g}{\e}+\d^{1/2}+\d^{1/4}+\e+(\frac{\g}{\d})^{1/2}\right)\left(1+|x_0|_\mH^2+|y_0|_\mH^2\right).
\de
By the above deduction, we infer that 
\ce
&&\mP\left(\sup\limits_{t\in[0,T]}\left|X^{\e,\g,u_{\e}}_t-\bar{X}^{u_\e}_{t}\right|_\mH>\eta\right)\\
&\leq&\frac{C}{R}\left(1+|x_0|_\mH^2+|y_0|_\mH^2\right)+ \frac{C_{R}}{\eta^2}\left(\frac{\g}{\e}+\d^{1/2}+\d^{1/4}+\e+(\frac{\g}{\d})^{1/2}\right)\left(1+|x_0|_\mH^2+|y_0|_\mH^2\right).
\de

Next, we take $\d=\g^{\iota}$ for $0<\iota<1$ and obtain that $\frac\g \e \rightarrow 0, \g\rightarrow 0, \d\rightarrow 0, \frac\g \d\rightarrow 0$ as $\e\rightarrow 0$, and further 
\ce
\lim _{\varepsilon \rightarrow 0}\mP\left(\sup\limits_{t\in[0,T]}\left|X^{\e,\g,u_{\e}}_t-\bar{X}^{u_\e}_{t}\right|_\mH>\eta\right)=0,
\de
as $\e\rightarrow0$ first and then $R\rightarrow\infty$. The proof is complete.
\end{proof}

Now, it is the position to prove Theorem \ref{fldpth}.

{\bf Proof of Theorem \ref{fldpth}.} 
By Proposition \ref{just1} and \ref{just2}, we know that Condition \ref{cond} holds. Thus, by Theorem \ref{ldpth}, we conclude that Theorem \ref{fldpth} is right.

\section{An example}\label{exam}

In this section, we give an example to explain our result.

\bx
Take $\mU_1=\mU_2=\mR$ and consider the following multiscale stochastic Burgers equations with reflection on the spatial domain $[0,1]$:
\be\left\{\begin{array}{ll}
\frac{\p X_t^{\e,\g}(\xi)}{\p t}=\frac{\p^2 X_t^{\e,\g}(\xi)}{\p\xi^2}+\frac{1}{2}\frac{\p(X_t^{\e,\g})^2(\xi)}{\p \xi}+f(X_t^{\e,\g},Y_t^{\e,\g})(\xi)+\sqrt \e g(X_t^{\e,\g})\frac{\p W^1(t,\xi)}{\p t}+\frac{\p K^{1,\e,\g}(t,\xi)}{\p t},\\
\frac{\p Y_t^{\e,\g}(\xi)}{\p t}=\frac{1}{\g}\left[\frac{\p^2 Y_t^{\e,\g}(\xi)}{\p\xi^2}+X_t^{\e,\g}(\xi)+l_1 Y_t^{\e,\g}(\xi)\right]+\frac{1}{\sqrt \g}l_2\frac{\p W^2(t,\xi)}{\p t}+\frac{\p K^{2,\e,\g}(t,\xi)}{\p t},\\
X_t^{\e,\g}(\xi)\geq 0, \quad Y_t^{\e,\g}(\xi)\geq 0, \\
X_t^{\e,\g}(0)=X_t^{\e,\g}(1)=Y_t^{\e,\g}(0)=Y_t^{\e,\g}(1)=0, \quad t>0,\\
X_0^{\e,\g}(\xi)=x_0(\xi)\geq 0, \quad Y_0^{\e,\g}(\xi)=y_0(\xi)\geq 0,
\end{array}
\right.
\label{exameq}
\ee
where $f: \mH\times\mH\to\mH, g: \mH\to\cL_2(\mR,\mH)$ are Lipschitz continuous and $l_1, l_2\in\mR, l_1, l_2\neq 0$, $\l_1-\frac{3}{2}|l_1|>0$. So, we observe that for $x,y\in\mH$,
\ce
&&F_1(x,y)=f(x,y), \quad G_1(x)=g(x),\\
&&F_2(x,y)=x+l_1y, \quad G_2(x,y)=l_2.
\de
Then, $({\bf H}_{F_1, G_1})$ and $({\bf H}^1_{F_2, G_2})$ are satisfied with $L_1=Lip(f), L_2=Lip(g), L_3=1, L_4=|l_1|, L_5=L_6=0$, where $Lip(f)$ denotes the Lipschitz constant for $f$. Note that $L_7=|l_2|>0$ and $\l_1-\frac{3}{2}L_4-L_6^2=\l_1-\frac{3}{2}|l_1|>0$. Thus, $({\bf H}^2_{G_2})$ and $({\bf H}^3_{F_2, G_2})$ hold and Theorem \ref{wellpose} implies that the system (\ref{exameq}) has a unique solution $(X^{\e,\g}, K^{1,\e,\g},Y^{\e,\g}, K^{2,\e,\g})$. 

Next, by $L_8=L_7=|l_2|>0$, we know that $({\bf H}^4_{G_2})$ is right. In terms of Theorem \ref{fldpth}, we infer that if $\lim\limits_{\e\to0}\g/\e=0$, $\{X^{\e,\g}, \e>0\}$ satisfies the LDP in $C([0,T], \mH)$ with the rate function $\Lambda$ given by
\ce
\Lambda(\phi)=\inf\limits _{\substack{\{h\in L^2([0,T],\sU): \\\phi=\bar{X}^{h}\}}}\left\{\frac{1}{2}\int_0^T|h(t)|_{\sU}^2\dif t\right\}, \quad \phi\in C([0,T], \mH),
\de
where $(\bar{X}^{h},\bar K^{1,h})$ solves the following stochastic Burgers equations with reflection 
\ce\left\{\begin{array}{ll}
\dif \bar X_t^{h}=[A\bar X_t^{h}+B(\bar X_t^{h},\bar X_t^{h})+\bar f(\bar X_t^{h})]\dif t+g(\bar X_t^{h})\Pi_1h(t)\dif t+\dif \bar K^{1,h}_t,\\
\bar X_t^{h}(\xi)\geq 0, \quad \xi\in[0,1],\\
\bar X_t^{h}(0)=\bar X_t^{h}(1)=0,\\
\bar X_0^{h}=x_0,
\end{array}
\right.
\de
$\bar f(x)=\int_\mH f(x,y)\mu^x(\dif y)$ and $\mu^x$ is the unique invariant measure of the following SPDE with reflection:
\ce\left\{\begin{array}{ll}
\dif Y_t^{x,y_0}=[AY_t^{x,y_0}+x+l_1Y_t^{x,y_0}]\dif t+l_2\dif \tilde W^2_t+\dif K^{2,x,y_0}_t,\\
Y_t^{x,y_0}(\xi)\geq 0,\quad \xi\in[0,1],\\
Y_t^{x,y_0}(0)=Y_t^{x,y_0}(1)=0,\\
Y_0^{x,y_0}=y_0.
\end{array}
\right.
\de

Finally, we mention that if $K^{1,\e,\g}=K^{2,\e,\g}=0$, the system (\ref{exameq}) falls into the framework in \cite{swxy}. Moreover, our LDP result is the same to the conclusion there.
\ex

\section{Appendix}\label{app}

In this section, we collect the properties of $B$ and $b$ (cf. \cite{zbz}).

\bl\label{Bbprop1}
$(i)$ For all $x, y, z \in \mV$,
$$
{_{\mV^*}}\<B(x, y), z\>_\mV=b(x,y,z)=-b(x,z,y)=-{_{\mV^*}}\< B(x, z), y\>_\mV.
$$

$(ii)$ For all $x, y, z \in \mV$,
$$
|{_{\mV^*}}\<B(x,y), z\>_\mV|=|b(x, y, z)| \leq 2\|x\|_\mV^{\frac{1}{2}}|x|_\mH^{\frac{1}{2}}\|z\|_{\mV}^{\frac{1}{2}}|z|_\mH^{\frac{1}{2}}\|y\|_\mV.
$$
\el

By the above lemma, we can obtain the following corollary.

\bc\label{Bbprop2}
\ce
&&b(x, y, y)=0, \quad {_{\mV^*}}\< B(x, y),y\>_\mV=0, \quad x, y \in \mV, \\
&&\|B(x, x)\|_{V^*} \leq 2\|x\|_\mV|x|_\mH, \quad x \in \mV.
\de
\ec

\bl\label{Bbprop3}
For any $x,y\in\mV$, it holds that
\ce
\|B(x,x)-B(y,y)\|_{\mV^*}\leq 2|x-y|_\mH(\|x\|_\mV+\|y\|_\mV).
\de
\el


\begin{thebibliography}{999}

\bibitem{btwy} L. Bo, D. Tang, Y. Wang and X. Yang: On the conditional default probability in a regulated market: a structural approach, {\it Quant. Finance}, 11(2011)1695-1702.

\bibitem{blz} Z. Brze\'zniak, Q. Li and T. Zhang: Large deviation principle of stochastic evolution equations with reflection, {\it J. Evol. Equ.}, 24(2024)91.

\bibitem{bm} Z. Brze\'zniak and E. Motyl: Existence of a martingale solution of the stochastic Navier-Stokes equations in unbounded 2D and 3D domains, {\it J. Differ. Equ.}, 254(2013)1627-1685.

\bibitem{zbz} Z. Brze\'zniak and T. Zhang: Reflection of stochastic evolution equations in infinite dimensional domains, {\it Ann. Inst. Henri Poincar\'e Probab. Stat.}, 59(2023)1549-1571.

\bibitem{abz} A. Budhiraja and P. Zoubouloglou: Large deviations for small noise diffusions over long time, {\it Transactions of the American Mathematical Society, Series B}, 11(2024)1-63.

\bibitem{css} Y. Chen, Y. Shi and X. Sun: Averaging principle for slow-fast stochastic Burgers equation driven by stable process, {\it Applied Mathematics Letters}, 103(2020)106199.

\bibitem{dsxz} Z. Dong, X. Sun, H. Xiao and J. Zhai: Averaging principle for one dimensional stochastic Burgers equation, {\it J. Differ. Equ.}, 265(2018)4749-4797.

\bibitem{de} P. Dupuis and R. Ellis: {\it A Weak Convergence Approach to the Theory of Large Deviations}, Wiley, New York, 1997.

\bibitem{ffk} J. Feng, J. Fouque and R. Kumar: Small-time asymptotics for fast mean-reverting stochastic volatility models, {\it Ann. Appl. Probab.}, 22(2012)1541-1575.

\bibitem{fg} F. Flandoli and D. Gatarek: Martingale and stationary solution for stochastic Navier-Stokes equations, {\it Probab. Theory Relat. Fields}, 102(1995)367-391.

\bibitem{fo} T. Funaki and S. Olla: Fluctuations for $\triangledown\phi$ interface model on a wall, {\it Stoch. Process. Appl.}, 94(2001)1-27.

\bibitem{ghl} J. Gao, W. Hong and W. Liu: Small noise asymptotics of multi-scale McKean-Vlasov stochastic dynamical systems, {\it J. Differ. Equ.}, 364(2023)521-575.

\bibitem{hkn} B. Hambly, J. Kalsi and J. Newbury: Limit order books, diffusion approximations and reflected SPDEs: From microscopic to macroscopic models, {\it Applied Mathematical Finance}, 27(2020)132-170.

\bibitem{hp} U. Haussmann and E. Pardoux: Stochastic variational inequalities of parabolic type, {\it Appl. Math. Optim.}, 20(1989)163-192.

\bibitem{hqs} M. Hong, Z. Qiu and C. Sun: Large deviations for the two-time scale 2D stochastic electrokinetic flow, {\it Potential Analysis},  63(2025)473-511.

\bibitem{hlll} W. Hong, M. Li, S. Li and W. Liu: Large deviations and averaging for stochastic tamed 3D Navier-Stokes equations with fast oscillations, {\it Appl. Math. Optim.}, 86(2022)15.

\bibitem{hll} W. Hong, S. Li and W. Liu: Freidlin-Wentzell type large deviation principle for multiscale locally monotone SPDEs, {\it SIAM J. Math. Anal.}, 53(2021)6517-6561.

\bibitem{hlls} W. Hong, S. Li, W. Liu, X. Sun: Central limit type theorem and large deviations for multi-scale McKean-Vlasov SDEs, {\it Probab. Theory Relat. Fields}, 187(2023)133-201.

\bibitem{hss} W. Hu, M. Salins and K. Spiliopoulos: Large deviations and averaging for systems of slow-fast stochastic reaction-diffusion equations, {\it Stoch. PDE Anal. Comp.}, 7(2019)808-874.

\bibitem{kk} H. Kang and T. Kurtz: Separation of time-scales and model reduction for stochastic reaction networks, {\it Ann. Appl. Probab.}, 23(2010)164-187.

\bibitem{yK} Y. Kifer: Large Deviations and Adiabatic Transitions for Dynamical Systems and Markov Processes in Fully Coupled Averaging, {\it Memoirs of the Amer. Math. Soc.}, 944, AMS, Providence, RI, 2009.

\bibitem{kp} R. Kumar and L. Popovic: Large deviations for multi-scale jump-diffusion processes, {\it Stoch. Process. Appl.}, 127(2017)1297-1320.

\bibitem{ku} H. J. Kushner: Large deviations for two-time-scale diffusions with delays, {\it Appl. Math. Optim.}, 62(2010)295-322.

\bibitem{mll} M. Li and W. Liu: Large deviation principle for multi-scale stochastic systems with monotone coefficients, {\it Communications in Mathematics and Statistics}, 14(2026)247-283.

\bibitem{tll} T. Li and F. Lin: Large deviations for two-scale chemical kinetic processe, {\it Commun. Math. Sci.}, 15(2017)123-163.

\bibitem{rL} R. Liptser: Large deviations for two scaled diffusions, {\it Probab. Theory Relat. Fields}, 106(1996)71-104.

\bibitem{jll} J. Liu and G. Long: Large deviations principle for multi-scale SDEs in H\"older norm, {\it Discrete Contin. Dyn. Syst. Ser. B}, 34(2026)141-166.

\bibitem{lrsx} W. Liu, M. R\"ockner, X. Sun and Y. Xie: Strong averaging principle for slow-fast stochastic partial differential equations with locally monotone coefficients, {\it Appl. Math. Optim.},  87(2023)39.

\bibitem{my} Z. Ma and J. Yang: Averaging principle for two time-scales stochastic partial differential equations with reflection, {\it Appl. Math. Optim.},  89(2024)59.

\bibitem{aw} A. Matoussi, W. Sabbagh and T. Zhang: Large deviation principles of obstacle problems for quasilinear stochastic PDEs, {\it Appl. Math. Optim.}, 83(2021)849-879.

\bibitem{ntt} P. Nguyen, K. Tawri and R. Temam: Nonlinear stochastic parabolic partial differential equations with a monotone operator of the Ladyzenskaya-Smagorinsky type, driven by a L\'evy noise, {\it J. Funct.
Anal.}, 281(2021)109157.

\bibitem{aP} A. Pilipenko: {\it An introduction to stochastic differential equations with reflection}. In: Lectures in Pure and Applied Mathematics, vol. 1, Potsdam University Press, Potsdam, 2014. 

\bibitem{aap} A. Puhalskii: On large deviations of coupled diffusions with time scale separation, {\it Ann. Probab.}, 44(2016)3111-3186.

\bibitem{q0} H. Qiao: Asymptotic behaviors of multiscale multivalued stochastic systems with small noises, to appear in {\it Proceedings of the Royal Society of Edinburgh Section A: Mathematics}, 2026.

\bibitem{q1} H. Qiao: Large deviations of multiscale multivalued McKean-Vlasov stochastic systems, {\it Stochastic Analysis and Applications}, 44(2026)565-604.

\bibitem{q2} H. Qiao: Large deviation principles for fully coupled multiscale multivalued stochastic systems, http://arxiv.org/abs/2512.10311.

\bibitem{q} H. Qiao: Uniform large deviation principles and averaging principles for stochastic Burgers type equations with reflection, {\it Nonlinear Analysis}, 274(2027)114263.

\bibitem{aS} A. Skorokhod: Stochastic equations for diffusion processes in a bounded region 1, 2, {\it Theory Probab. Appl.}, 6(1961)264-274, 7(1962)3-23.

\bibitem{ks1} K. Spiliopoulos: Large deviations and importance sampling for systems of slow-fast motion, {\it Appl. Math. Optim.}, 67(2013)123-161.

\bibitem{swxy} X. Sun, R. Wang, L. Xu and X. Yang: Large deviation for two-time-scale stochastic Burgers equation, {\it Stochastics and Dynamics}, 21(2021)2150023(37 pages).

\bibitem{aV1} A. Veretennikov: On large deviations for SDEs with small diffusion and averaging, {\it Stoch. Process. Appl.}, 89(2000)69-79.

\bibitem{aV2} A. Veretennikov: On large deviations in the averaging principle for SDEs with a ``full dependence" revisited, {\it Discrete Contin. Dyn. Syst. Ser. B}, 18(2013)523-549.

\bibitem{wzz} R. Wang, J. Zhai and S. Zhang: Large deviation principle for stochastic Burgers type equation with reflection, {\it Communications on Pure and Applied Analysis}, 21(2022)213-238.

\bibitem{wrd} W. Wang, A. J. Roberts and J. Duan: Large deviations and approximations for slow-fast stochastic reaction-diffusion equations, {\it J. Differ. Equ.}, 253(2012)3501-3522.

\bibitem{wtry} F. Wu, T. Tian, J. Rawlings and G. Yin: Approximate method for stochastic chemical kinetics with two-time scales by chemical Langevin equations, {\it J. Chem. Phys.}, 144(2016)174112.

\bibitem{xz} T. Xu and T. Zhang: White noise driven SPDEs with reflection: existence, uniqueness and large deviation principles, {\it Stoch. Process. Appl.}, 119(2009)3453-3470.

\bibitem{xxyp} W. Xu, Y. Xu, X. Yang and B. Pei: Large deviation principle for slow-fast systems with infinite-dimensional mixed fractional Brownian motion, {\it Journal of Statistical Physics}, 192(2025)161.

\bibitem{yz} J. Yang and T. Zhang: Existence and uniqueness of invariant measures for SPDEs with two reflecting walls, {\it J. Theoret. Probab.}, 27(2014)863-877.

\bibitem{ytw} X. Yin, Y. Tian and J.-L. Wu: Large deviation and averaging for multi-scale stochastic tidal dynamics equation, {\it J. Math. Phys.}, 66(2025)072702.

\bibitem{zhangt0} T. Zhang: Lattice approximations of reflected stochastic partial differential equations driven by space-time white noise, {\it Ann. Appl. Probab.}, 26(2016)3602-3629.

\bibitem{zhangt} T. Zhang: Stochastic Burgers type equations with reflection: existence, uniqueness, {\it J. Differ. Equ.}, 267(2019)4537-4571.

\end{thebibliography}
\end{document}